\documentclass[a4paper,10pt,reqno]{elsarticle}

\usepackage{amsmath, amstext, paralist, amsthm, amssymb, epsfig, paralist, tikz, tikz-cd, graphicx}
\usepackage{hyperref}
\usepackage{xurl}
\usepackage{float}

\newtheorem{thm}{Theorem}[section]
\newtheorem{lem}[thm]{Lemma}
\newtheorem{prop}[thm]{Proposition}

\newtheorem{cor}[thm]{Corollary}

\newcommand{\alt}{\operatorname{alt}\nolimits}
\newcommand{\Aut}{\operatorname{Aut}\nolimits}
\newcommand{\emb}{\operatorname{emb}\nolimits}
\newcommand{\out}{\operatorname{out}\nolimits}
\begin{document}
	
	\begin{frontmatter}
		
		\title{Two Types of Cards in Kelly-Ulam Conjecture}
		
		\author{Max Zinchenko\fnref{myfootnote}}
		\affiliation{organization={Novosibirsk State University},
			addressline={1 Pirogova St.},
			city={Novosibirsk},
			postcode={630090},
			country={Russia}}
		
		\cortext[mycor]{Corresponding author}
		\fntext[myfootnote]{Email: m.zinchenko@g.nsu.ru}
		
		\begin{abstract}
			We consider the case in Kelly-Ulam conjecture where there are only two types of cards in deck. It was possible to prove that graphs with this property are uniquely reconstructed by their card deck.
		\end{abstract}
		
		
	\end{frontmatter}
	
\section{Introduction}
The Graph Reconstruction Conjecture is one of the most famous and robust open problems in graph theory. The conjecture asserts that any finite simple graph with at least three vertices is uniquely determined (up to isomorphism) by its multiset of vertex-deleted subgraphs, which is called a ``deck.''
\par The problem was formulated in Paul Kelly's dissertation \cite{Kelly1942}. Stanisław Ulam included this problem in his book \cite{Ulam1960}; therefore, it is often called the Kelly--Ulam conjecture. Although this problem is still open, the conjecture has been proven for many important graph classes, for instance:
\begin{enumerate}
	\item Trees \cite{Kelly1957}
	\item Disconnected graphs \cite{Harary1964}
	\item Regular graphs \cite{BondyHemminger1977}
\end{enumerate}
Furthermore, computer checks, such as those by Brendan McKay \cite{McKay2022}, have shown that the conjecture is true for all graphs with up to 13 vertices.
\par Initially, the idea was to describe all graphs with reconstruction number equal to $3$, because according to \cite{Bollobas1990}, almost every graph is reconstructed by just 3 cards.
\par Let $t(G)$ denote the number of isomorphism classes among the cards of $G$. In the process of studying graphs with $r(G)=3$, I formulated the conjecture that $t(G) \geq 3 \Longrightarrow r(G)=3$, which we have been unable to refute by computer brute force and have failed to prove. Therefore, I introduce a new method by which we can attempt to prove the Kelly--Ulam conjecture:
\begin{enumerate}
	\item If $t(G) = 1$, then $G$ is regular, and this case has been proven;
	\item The case $t(G)=2$ is proved in this paper;
	\item Prove the conjecture $t(G) \geq 3 \Longrightarrow r(G)=3$.
\end{enumerate}
\subsection{Organization of the paper}
The main result of this paper is the following theorem:
\begin{thm}
	Let $G$ be a finite, simple, undirected graph such that $|G| \geq 3$ and $t(G)=2$. Then $G$ is uniquely determined by its deck.
\end{thm} 

The remainder of this paper is organized as follows. 
\par In Section $2$, we perform a reduction of the general case. Assuming that a counterexample exists, we show that it must be a graph with only two vertex degrees, $d$ and $d+1$. From this point onward, we treat this simplified case as our main setting.
\par In Section $3$, we introduce the notion of an extension of a card in the deck.
\par In Section $4$, we show that in our main case, the graph has exactly two extensions.
\par In Sections $5$--$9$, we define functions and algebraic invariants associated with these extensions and prove several identities.
\par In Section $10$, we demonstrate that one of these identities is impossible, which implies that a counterexample does not exist.

\subsection{Preliminaries}
All graphs are finite, simple, and undirected. For a graph $G$, let $V(G)$ be its vertex set and put $n=|V(G)|$. Adjacency in $G$ is denoted by $x\sim_G y$, and the neighborhood of a vertex $x$ is denoted by $N_G(x)$.

The card $G-v$ is the subgraph induced by $V(G)\setminus\{v\}$. The deck of $G$ is the multiset of all its cards. A graph $H$ on the same vertex set $V$ has the same deck as $G$ (written $G\sim H$) if there is a bijection $\sigma\colon V\to V$ such that
$$
G-v \cong H-\sigma(v)\quad\text{for every }v\in V.
$$

For a graph $G$, let $A_G$ be its adjacency matrix and let $a_G(u)$ be the adjacency vector of a vertex $u$ in $G$. For functions $f,g\colon V\to\mathbb{Q}$, define
$$
\langle f,g\rangle := \sum_{x\in V} f(x)g(x).
$$

Let $\varepsilon_u$ be the indicator of the vertex $u$; in other words,
$$
\varepsilon_u(x)=
\begin{cases}
	1, & \text{if } x=u,\\
	0, & \text{if } x\neq u.
\end{cases}
$$
The notation $[P]$ equals $1$ if $P$ is true, and $0$ otherwise. Finally, $\mathbf{1}$ denotes the all-ones vector.
	\section{Reducing to main case}
\begin{lem}\label{l2.1}
Let $n \ge 3$ and let $G \sim H$ be witnessed by a bijection $\sigma$. Then $|E(G)| = |E(H)|$ and $\deg_H(\sigma(v)) = \deg_G(v)$ for all $v$. Furthermore, if $G - u \cong G - u'$, then $\deg(u) = \deg(u')$.

\begin{proof}
	Every edge is contained in exactly $n - 2$ cards, so
	\[
	\sum_{v \in V} |E(G - v)| = (n - 2)|E(G)|,
	\]
	and similarly for $H$. Since $|E(G - v)| = |E(H - \sigma(v))|$ for all $v$, and $n - 2 > 0$, we deduce that $|E(G)| = |E(H)|$.
	
	Further, for every $v$,
	\[
	\deg_G(v) = |E(G)| - |E(G - v)| = |E(H)| - |E(H - \sigma(v))| = \deg_H(\sigma(v)).
	\]
	
	The last assertion follows directly from the same identity: if $G - u \cong G - u'$, then $|E(G - u)| = |E(G - u')|$, and hence
	\[
	\deg(u) = |E(G)| - |E(G - u)| = |E(G)| - |E(G - u')| = \deg(u').
	\]
\end{proof}

\end{lem}
	
\begin{lem}\label{l2.2}
	Let $n \ge 3$ and suppose there exists a vertex $v$ with the following property: if $x, y \neq v$ have the same degree in the card $G - v$, then either both are adjacent to $v$ or both are non-adjacent to $v$. Then every graph $H \sim G$ is isomorphic to $G$.
	
\begin{proof}
	Let $\sigma$ be the bijection from the definition of $G \sim H$. Then there exists a bijection $\sigma$ with $v' = \sigma(v)$, and an isomorphism $\varphi: G - v \to H - v'$. For $x \neq v$, define
	\[
	c(x) = \deg_{G-v} x = \deg_{H-v'} \varphi(x), \quad
	s(x) = \begin{cases} 1, & \text{if } x \sim_G v, \\ 0, & \text{if } x \nsim_G v, \end{cases} \quad
	t(x) = \begin{cases} 1, & \text{if } \varphi(x) \sim_H v', \\ 0, & \text{otherwise.} \end{cases}
	\]
	Then
	\[
	\deg_G x = c(x) + s(x), \qquad \deg_H \varphi(x) = c(x) + t(x).
	\]
	By Lemma \ref{l2.1}, the degree multisets of $G$ and $H$ coincide, and $\deg_H v' = \deg_G v$. Since $\varphi$ is a bijection from $V \setminus \{v\}$ to $V \setminus \{v'\}$, for every $j$, the number of $x \neq v$ such that $c(x) + s(x) = j$ equals the number of $x \neq v$ with $c(x) + t(x) = j$ (both are equal to the number of vertices of degree $j$ minus $[\deg v = j]$).
	
	We prove by induction on $k$ that
	$$
	\#\{x: c(x) = k, s(x) = 1\} = \#\{x: c(x) = k, t(x) = 1\}.
	$$
	We have
	$$
	\#\{c + s = j\} = \#\{c = j, s = 0\} + \#\{c = j - 1, s = 1\},
	$$
	$$
	\#\{c = j, s = 0\} + \#\{c = j, s = 1\} = \#\{c = j\},
	$$
	and analogously for $t$.
	
	For $j = 0$, the second term $\#\{c = j - 1, s = 1\}$ is zero, hence
	$$
	\#\{c = 0, s = 0\} = \#\{c = 0, t = 0\},
	$$
	which implies
	$$
	\#\{c = 0, s = 1\} = \#\{c = 0, t = 1\}.
	$$
	For the inductive step $k \to k+1$, assume
	$$
	\#\{c = k, s = 1\} = \#\{c = k, t = 1\}.
	$$
	Then
	$$
	\#\{c = k+1, s = 0\} = \#\{c = k+1, t = 0\}.
	$$
	By the lemma's conditions, $s$ is constant on each class $C_k = \{x: c(x) = k\}$ (since it cannot be that $y \in C_k$ is adjacent to $v$ while $x \in C_k$ is not).
	
	If $s \equiv 1$ on $C_k$, then $\#\{x \in C_k : t(x) = 1\} = |C_k|$, i.e., $t \equiv 1$ on $C_k$; if $s \equiv 0$, then $t \equiv 0$. Thus $t \equiv s$, which means
	$$
	x \sim_G v \iff \varphi(x) \sim_H v'.
	$$
	Extending $\varphi$ by $v \mapsto v'$ yields an isomorphism $G \cong H$. 

\end{proof}
	
\end{lem}	

\begin{prop}[Reduction]\label{p2.3}
	Let $n \ge 3$ and $t(G) = 2$, and suppose there exists $H \sim G$ with $H \not\cong G$. Then for some $d$:
	\begin{enumerate}
		\item[(i)] the degree of every vertex is either $d$ or $d+1$, and both degrees occur;
		\item[(ii)] every vertex $v$ has a neighbor of degree $d+1$ and a non-neighbor $x \neq v$ of degree $d$.
	\end{enumerate}
	
	\begin{proof}
		(i) By Lemma \ref{l2.1}, vertices of different degrees have non-isomorphic cards. Since $t(G) = 2$, there are  two distinct degrees in $G$.
		
		Suppose, for the sake of contradiction, that no two distinct degrees differ by exactly $1$. We claim that the hypothesis of Lemma \ref{l2.2} holds for every vertex $v$. Indeed, the degree of a vertex $x \neq v$ in the card $G - v$ is given by
		$$
		\deg_{G-v}(x) = \deg_G(x) - s(x),
		$$
		where $s(x)$ from previous lemma. Thus, the conditions of Lemma \ref{l2.2} can be rewritten as:
		$$
		\text{if } \deg_G(x) - s(x) = \deg_G(y) - s(y), \text{ then } s(x) = s(y).
		$$
		
		Assume this condition fails. Then there exist $x, y \neq v$ such that
		$$
		\deg_G(x) - s(x) = \deg_G(y) - s(y), \quad \text{but} \quad s(x) \neq s(y).
		$$
		Without loss of generality, let $s(x) = 1$ and $s(y) = 0$. Then
		$$
		\deg_G(x) - 1 = \deg_G(y) \quad\Longrightarrow\quad \deg_G(x) - \deg_G(y) = 1,
		$$
		which contradicts our assumption that no two degrees differ by exactly $1$. Hence, the conditions of Lemma \ref{l2.2} is satisfied for every $v$.
		
		Applying Lemma \ref{l2.2}, we conclude that every graph $H \sim G$ is isomorphic to $G$. However, this contradicts the assumption of the proposition that there exists $H \sim G$ with $H \not\cong G$. Therefore, our supposition was false, and there must exist two distinct degrees differing by exactly $1$. Since there are at most two distinct degrees, they must be $d$ and $d+1$.
		\par (ii) Call vertices of degree $d$ \textit{low} and vertices of degree $d+1$ \textit{high}.
		
		Consider a fixed vertex $v$ and the card $G-v$. For every $x \neq v$,
		$$
		\deg_{G-v}(x)=\deg_G(x)-s(x).
		$$
		Which degrees occur in $G-v$?
		
		\begin{table}[H]
			\centering
			\begin{tabular}{lccc}
				\hline
				Type of vertex $x$ & $\deg_G(x)$ & $s(x)$ & $\deg_{G-v}(x)$ \\
				\hline
				Low, non-adjacent to $v$ & $d$ & $0$ & $d$ \\
				Low, adjacent to $v$ & $d$ & $1$ & $d-1$ \\
				High, non-adjacent to $v$ & $d+1$ & $0$ & $d+1$ \\
				High, adjacent to $v$ & $d+1$ & $1$ & $d$ \\
				\hline
			\end{tabular}
		\end{table}
		The possible degrees in $G-v$ are $d-1, d, d+1$, but degree $d$ is ambiguous. If $v$ has no high neighbours or is adjacent to all low vertices, then the degree in the card uniquely determines adjacency to $v$; hence for all $x,y \neq v$ with equal degrees in $G-v$, they are either both adjacent to $v$ or both non-adjacent to $v$. Thus the condition of Lemma \ref{l2.2} is satisfied, but this contradicts $H \not\cong G$.
	\end{proof}
	
\end{prop}
Further $L$ and $U$ - sets of low and high vetices, $p= |L| \ge 0, q=|U| \ge 0$.

\begin{lem}\label{l2.4}
	Suppose that $t(G)= 2$ and that condition~(i) holds. Then all cards $G-v$, $v\in L$, are pairwise isomorphic, and so are all cards $G-w$, $w\in U$. There exist numbers $b,e$ such that every low vertex has exactly $b$ high neighbours and every high vertex has exactly $e$ high neighbours. Every graph $K\sim G$ has degrees $d,d+1$ with the same multiplicities $p,q$, the same parameters $b,e$, and the same degree sum.
	
	\begin{proof}
		There are two degree classes and two isomorphism types of cards, and cards of vertices of different degrees are non-isomorphic; hence each degree class corresponds to exactly one card type. In the card $G-v$, the vertices of degree $d+1$ are precisely the high non-neighbours of $v$. Their number is $q-b(v)$ for $v\in L$ and $q-1-e(v)$ for $v\in U$, where $b(v)$ and $e(v)$ denote the numbers of high neighbours of $v$. These quantities are invariants of the card type, and therefore $b(v)$ and $e(v)$ are constant. For $K\sim G$, Lemma \ref{l2.1} gives $\deg_K\sigma(v)=\deg_G v$ and $K-\sigma(v)\cong G-v$; the same counting argument yields the same values $b,e$ for $K$. Finally, the degree sums are equal: $2|E(G)|=2|E(K)|$.
	\end{proof}
\end{lem}

We set
$$
\theta = e - b, \qquad
c = d + 1 - e \ \ \text{(the number of low neighbours of a high vertex)},
pb = qc.
$$
In particular, condition (ii) implies that $b \ge 1$.

A \emph{counterexample} is a pair of graphs $(G,H)$ on $V$ such that $n \ge 3$, $t(G) = 2$, conditions (i) and (ii) hold for some $d$, $G \sim H$, and $G \ncong H$.

\begin{lem}\label{l2.5}
	If $(G,H)$ is a counterexample with parameters $d,b,e$, then $(\overline{G},\overline{H})$ is a counterexample with $d' = n-2-d$; the low and high classes are interchanged, and $\theta' = -\theta$.
\end{lem}

\begin{proof}
	Since $\overline{G-v} = \overline{G}-v$, we have $\overline{G} \sim \overline{H}$, $t(\overline{G}) = t(G)$, and $\overline{G} \ncong \overline{H}$. The degrees in the complement are $n-1-\deg$: the high vertices of $G$ acquire degree $d'$, and the low ones acquire degree $d'+1$. A neighbour of degree $d'+1$ in $\overline{G}$ is a low non-neighbour in $G$; a non-neighbour $x \ne v$ of degree $d'$ in $\overline{G}$ is a high neighbour in $G$. Hence condition~(ii) for $\overline{G}$ coincides with condition~(ii) for $G$. Finally, $b' = p - c$ (the low non-neighbours of a high vertex), $e' = p - 1 - (d-b)$ (the low non-neighbours $x \ne v$ of a low vertex $v$), and
	$$
	\theta' = e' - b' = c + b - 1 - d = b - e = -\theta. \qedhere
	$$
\end{proof}

\section{Card extensions}

Let $\mathcal{K}$ be the set of isomorphism classes of graphs $K$ on $V$ such that $K \sim G$. By Lemma \ref{l2.4}, every $K \in \mathcal{K}$ has the same vertex degrees, the same parameters $b$ and $e$, and the same sum of degrees as $G$.

For $K \in \mathcal{K}$, define
\[
\chi_K(x) = \deg_K(x) - d;
\]
in other words, $\chi_K$ is the indicator of the set of high vertices in $K$. Let $V = \{v_1, \dots, v_n\}$ and put
\[
\chi_K =
\begin{pmatrix}
	\chi_K(v_1) \\
	\vdots \\
	\chi_K(v_n)
\end{pmatrix}.
\]
It is easy to see that
\begin{equation}
	A_K \chi_K = b \mathbf{1} + \theta \chi_K, \tag{1}
\end{equation}
where $\theta = e - b$ and \(\mathbf{1}\) denote the all-ones vector. Indeed,
\[
(A_K \chi_K)(x) = \sum_{y \sim_K x} \chi_K(y),
\]
which is exactly the number of high neighbours of the vertex $x$.

Now consider the right-hand side of~(1):
\begin{itemize}
	\item if $x$ is low, i.e.\ $\chi_K(x) = 0$, then the right-hand side equals $b$;
	\item if $x$ is high, i.e.\ $\chi_K(x) = 1$, then the right-hand side equals $e$.
\end{itemize}

Let $K \in \mathcal{K}$ and $u \in V$. An \textit{extension} of the card $K - u$ is a graph $K' \in \mathcal{K}$ that coincides with $K$ on all pairs not containing $u$, i.e for all $x,y \in V \setminus \{u\}$ we have $\{x,y\} \in E(K) \iff \{x,y\} \in E(K')$. Denote the set of extensions by $C_K(u)$. Clearly, $K \in C_K(u)$ and $C_{K'}(u) = C_K(u)$ for $K' \in C_K(u)$. For graphs $K_1, K_2$ that coincide without $u$, we set
$$
z=z[K_1, K_2; u] = a_{K_1}(u) - a_{K_2}(u) \in \{-1, 0, 1\}^V,
$$
$$
k[K_1, K_2; u] = \#\{x : z[K_1, K_2; u](x) = 1\}.
$$
We have \(z(u) = 0\), and \(K_1 = K_2 \iff z = 0\).
In other words
\begin{itemize}
	\item $z(x) = 1$ if $x$ is a neighbour of $u$ in $K_1$ but not in $K_2$;
	\item $z(x) = -1$ if $x$ is a neighbour of $u$ in $K_2$ but not in $K_1$;
	\item $z(x) = 0$ if the adjacency status of $x$ to $u$ is the same in both graphs.
\end{itemize}
The number $k[K_1, K_2; u]$ is simply the number of vertices $x$ for which $z(x) = 1$, i.e., the number of ``extra'' edges incident to $u$ in $K_1$ compared to $K_2$.

\begin{lem}\label{l3.1}
	Let $K_1, K_2 \in \mathcal{K}$ coincide outside $u$, $z = z[K_1, K_2; u]$, $k = k[K_1, K_2; u]$, and $s = \chi_{K_1}(u)$. Then
	\begin{enumerate}
		\item[(a)] $\deg_{K_1} u = \deg_{K_2} u$, $\sum_x z(x) = 0$ (so $-1$ also occurs $k$ times and $k[K_2, K_1; u] = k$), and $\chi_{K_1} - \chi_{K_2} = z$;
		\item[(b)] if $z(x) = 1$, then $x$ is high in $K_1$; if $z(x) = -1$, then $x$ is low in $K_1$;
		\item[(c)] $(A_{K_1} - (\theta - s)I)z = k \varepsilon_u$;
		\item[(d)] if $K_1 \neq K_2$, then $k \ge 1$.
	\end{enumerate}
\end{lem}

\begin{proof}
	$(a)$ For $x \ne u$, only the pair $\{x,u\}$ can differ, and hence
	$\deg_{K_1} x - \deg_{K_2} x = z(x)$. By Lemma \ref{l2.4} graphs $K_1$ and $K_2$ have same vertex degrees sum, hence 
	$$\sum_{x \in V} \deg_{K_1}(x)=\sum_{x \in V} \deg_{K_2}(x) \Longrightarrow \sum_{x \in V} z(x) = 0$$
	By $K_1 - u = K_2 - u$ we have $\deg_{K_1} (u) = \deg_{K_2} (u)$. By definition of $\chi_K$ we have
	
	$$(\chi_{K_1} - \chi_{K_2})(x)=\deg_{K_1}(x)-d-\deg_{K_2}+d = z(x)$$
	$(b)$ follows from~$(a)$, since
	$\chi \in \{0,1\}$.\\
	$(c)$
	Let $x \neq u$. From equation~(1) applied to \(K_1\), for any vertex \(x\),
	\[
	(A_{K_1} \chi_{K_1})(x) = b + \theta \chi_{K_1}(x). \tag{2}
	\]
	Since \(z = \chi_{K_1} - \chi_{K_2}\), we can write \(\chi_{K_1} = \chi_{K_2} + z\). Substituting this into the left-hand side of~(2) gives
	\[
	(A_{K_1} (\chi_{K_2} + z))(x) = b + \theta \bigl(\chi_{K_2}(x) + z(x)\bigr).
	\]
	Expanding the left-hand side, we get
	\[
	(A_{K_1} \chi_{K_2})(x) + (A_{K_1} z)(x) = b + \theta \chi_{K_2}(x) + \theta z(x). \tag{3}
	\]
	
	Now we need to express \((A_{K_1} \chi_{K_2})(x)\) in terms of \((A_{K_2} \chi_{K_2})(x)\). Since \(K_1\) and \(K_2\) coincide outside \(u\), their adjacency matrices \(A_{K_1}\) and \(A_{K_2}\) differ only in the column corresponding to \(u\). For any \(x \neq u\), the \(x\)-th row of \(A_{K_1}\) differs from the \(x\)-th row of \(A_{K_2}\) only at the entry in column \(u\). Specifically,
	\[
	A_{K_1}(x, u) - A_{K_2}(x, u) = z(x).
	\]
	where \(A_{K_i}(x, u)\) is an entry of the adjacency matrix of the graph \(K_i\), located at the intersection of the row corresponding to vertex \(x\) and the column corresponding to vertex \(u\).
	\[
	(A_{K_1} \chi_{K_2})(x) = \sum_{y \in V} A_{K_1}(x, y) \chi_{K_2}(y)
	= \sum_{y \neq u} A_{K_2}(x, y) \chi_{K_2}(y) + A_{K_1}(x, u) \chi_{K_2}(u).
	\]
	Since \(A_{K_1}(x, u) = A_{K_2}(x, u) + z(x)\), this becomes
	\[
	(A_{K_1} \chi_{K_2})(x) = \sum_{y \neq u} A_{K_2}(x, y) \chi_{K_2}(y) + \bigl(A_{K_2}(x, u) + z(x)\bigr) \chi_{K_2}(u).
	\]
	Adding and subtracting the missing term \(A_{K_2}(x, u) \chi_{K_2}(u)\) carefully, we obtain
	\[
	(A_{K_1} \chi_{K_2})(x) = (A_{K_2} \chi_{K_2})(x) + z(x) \chi_{K_2}(u). \tag{4}
	\]
	By~(1) applied to \(K_2\), we have \((A_{K_2} \chi_{K_2})(x) = b + \theta \chi_{K_2}(x)\). Also, by part~(a) of the lemma, \(\chi_{K_2}(u) = \chi_{K_1}(u) = s\). Substituting this into~(4) yields
	\[
	(A_{K_1} \chi_{K_2})(x) = b + \theta \chi_{K_2}(x) + s z(x). \tag{5}
	\]
	
	Now substitute~(5) into the left-hand side of~(3):
	\[
	\bigl[b + \theta \chi_{K_2}(x) + s z(x)\bigr] + (A_{K_1} z)(x) = b + \theta \chi_{K_2}(x) + \theta z(x).
	\]
	Cancelling the common terms \(b\) and \(\theta \chi_{K_2}(x)\) on both sides, we are left with
	\[
	s z(x) + (A_{K_1} z)(x) = \theta z(x).
	\]
	Rearranging the terms gives
	\[
	(A_{K_1} z)(x) = (\theta - s) z(x).
	\]
	This can be rewritten as
	\[
	\bigl(A_{K_1} - (\theta - s)I\bigr) z (x) = 0 \qquad \text{for all } x \neq u. \tag{6}
	\]

When \(x = u\)
	
	\[
	(A_{K_1} z)(u) = \sum_{y \in V} A_{K_1}(u, y) z(y)
	= \sum_{y \in V} a_{K_1}(u)(y) \bigl(a_{K_1}(u)(y) - a_{K_2}(u)(y)\bigr).
	\]
	Since \(a_{K_1}(u)(y)\) and \(a_{K_2}(u)(y)\) are either \(0\) or \(1\), the product \(a_{K_1}(u)(y)\bigl(a_{K_1}(u)(y) - a_{K_2}(u)(y)\bigr)\) equals \(1\) if and only if \(a_{K_1}(u)(y) = 1\) and \(a_{K_2}(u)(y) = 0\). In all other cases, it equals \(0\). Therefore, the sum counts exactly the number of vertices \(y\) that are neighbours of \(u\) in \(K_1\) but not in \(K_2\). By definition, this number is exactly \(k = k[K_1, K_2; u]\). Hence,
	\[
	(A_{K_1} z)(u) = k. \tag{7}
	\]
	Also, since \(z(u) = 0\) (as \(u\) has no loops and hence the \(u\)-th coordinate of \(z\) is zero), we have
	\[
	\bigl(A_{K_1} - (\theta - s)I\bigr) z (u) = (A_{K_1} z)(u) - (\theta - s) z(u) = k - (\theta - s) \cdot 0 = k. \tag{8}
	\]
	
	Recall that \(\varepsilon_u\) is the indicator vector of the vertex \(u\), i.e., \(\varepsilon_u(u) = 1\) and \(\varepsilon_u(x) = 0\) for \(x \neq u\). From~(6) and~(8), we conclude that for every vertex \(x \in V\),
	\[
	\bigl(A_{K_1} - (\theta - s)I\bigr) z (x) =
	\begin{cases}
		0 = k \cdot 0 = k \varepsilon_u(x), & \text{if } x \neq u, \\
		k = k \cdot 1 = k \varepsilon_u(u), & \text{if } x = u.
	\end{cases}
	\]
	In other words,
	\[
	\bigl(A_{K_1} - (\theta - s)I\bigr) z = k \varepsilon_u,
	\]
	$(d)$ If $k=0$, then by $(a)$ we have $z\equiv0$, hence $K_1 = K_2$. 	
\end{proof}

\begin{lem}\label{l3.2}
	For every $K \in \mathcal{K}$ and every $u \in V$, the class $\mathcal{C}_K(u)$ contains a graph isomorphic to $G$ and a graph isomorphic to $H$. In particular, there exists $K' \in \mathcal{C}_K(u)$ with $K' \neq K$.
\end{lem}

\begin{proof}
	Let $X$ is arbitrary graph on $V$ such that $\psi\colon K-u \to X-x$ be an isomorphism of cards. We define a graph $K_X$ on $V$ as follows: it coincides with $K$ outside $u$, and $u \sim y \iff \psi(y) \sim_X x$. The extension of $\psi$ by the map $u \mapsto x$ is an isomorphism $K_X \cong X$.
	
	Since $K \sim G \sim H$, the card $K-u$ is isomorphic to some card $G-v$ and to some card $H-v'$. Indeed $K_G \cong G$ and $K_H \cong H$; both graphs belong to $\mathcal{K}$ and coincide with $K$ outside $u$, that is, they belong to $\mathcal{C}_K(u)$. Since $G \ncong H$, we have $K_G \ncong K_H$, and hence at least one of them is distinct from $K$.
\end{proof}

\begin{lem}\label{l3.3}
	Let $K \in \mathcal{K}$, $u \in V$, $s = \chi_K(u)$, and suppose that $A_K f = (\theta - s) f$. Then $f(u) = 0$. Consequently, $f$ vanishes on the entire class $\{x : \chi_K(x) = s\}$.
\end{lem}

\begin{proof}
	Take $K' \neq K$ as in Lemma \ref{l3.2}, and let $z = z[K,K';u]$ and $k \ge 1$. By Lemma \ref{l3.1} and the symmetry of $A_K$, we obtain
	$$
	f^{\mathsf{T}}k \varepsilon_u=k f(u) = f^{\mathsf{T}} \bigl(A_K - (\theta - s) I\bigr) z
	= \Bigl(\bigl(A_K - (\theta - s) I\bigr) f\Bigr)^{\mathsf{T}} z = 0.
	$$
	The second assertion follows by applying the first one to every vertex of the class~$s$.
\end{proof}

\begin{lem}\label{l3.4}
	Let $K \in \mathcal{K}$, $u \in V$, $s = \chi_K(u)$, $K_1, K_2 \in \mathcal{C}_K(u) \setminus \{K\}$, $z_i = z[K,K_i;u]$, and $k_i = k[K,K_i;u]$. Then $k_1 = k_2$ and $z_1(x) = z_2(x)$ for all $x$ with $\chi_K(x) = s$.
\end{lem}

\begin{proof}
	By Lemma~\ref{l3.1}(c), the vector $y = k_2 z_1 - k_1 z_2$ satisfies
	$$
	\bigl(A_K - (\theta - s) I\bigr) y = k_2 k_1 \varepsilon_u - k_1 k_2 \varepsilon_u = 0;
	$$
	hence, by Lemma~\ref{l3.3}, $y = 0$ on the class~$s$. This class contains a vertex $x_0$ with $z_1(x_0) \ne 0$. 
	
	For $s = 0$, the vertex $u$ is low. Since $\sum_x z_1(x) = 0$ and $k_1 \geq 1$ by Lemma~\ref{l3.1}(a)(d), there exists a vertex $x_0$ such that $z_1(x_0) = -1$. By Lemma~\ref{l3.1}(b), $x_0$ is a low vertex, hence $\chi_{K}(x_0) = \chi_{K}(u) = s$.
	
	For $s=1$, analogously, there exists a vertex $x_1$ such that $z_1(x_1) = 1$, and this vertex is high, i.e., $\chi_K(x_1) = 1$.
	
	Since $k_2 z_1(x_0) = k_1 z_2(x_0)$, $z_i(x_0) \in \{-1, 0, 1\}$, and $k_i > 0$, it follows that $z_2(x_0) = z_1(x_0)$ and $k_1 = k_2$. Therefore $k_1(z_1 - z_2) = y = 0$ on the class~$s$.
\end{proof}

\begin{cor}\label{c3.5} \medspace
	\begin{enumerate}
		\item[(a)] For every $u$, the number $k[X,Y;u]$ is the same for all pairs $X \ncong Y$ in $\mathcal{C}_G(u)$; we denote it by $\kappa(u) \ge 1$.
		\item[(b)] If $\deg u = \deg u'$, then there exists a permutation $\pi$ of the set $V$ with $\pi(u) = u'$ such that the map $X \mapsto X^{\pi}$ maps $\mathcal{C}_G(u)$ bijectively onto $\mathcal{C}_G(u')$; in particular, $\kappa(u) = \kappa(u')$.
	\end{enumerate}
\end{cor}

\begin{proof}
	(a) For $X \in \mathcal{C}_G(u)$ we have $\mathcal{C}_X(u) = \mathcal{C}_G(u)$, and by Lemma~\ref{l3.4} (with $K = X$) the number $k[X,Y;u]$ does not depend on $Y \ncong X$. Moreover, $k[X,Y;u] = k[Y,X;u]$ (Lemma~\ref{l3.1}(a)). Any two pairs can be joined by a chain of pairs in which any two consecutive pairs share a common element. By Lemma~\ref{l3.2}, $|\mathcal{C}_G(u)| \ge 2$, and $\kappa(u) \ge 1$ by Lemma~\ref{l3.1}(d).
	
	(b) The cards $G-u$ and $G-u'$ are isomorphic (Lemma~\ref{l2.4}); this isomorphism extends to a permutation $\pi$ with $\pi(u) = u'$. If $X \in \mathcal{C}_G(u)$, then $X^{\pi} \cong X$ lies in $\mathcal{K}$, and for $x, y \ne u'$ we have
	$$
	x \sim_{X^{\pi}} y \iff \pi^{-1}(x) \sim_X \pi^{-1}(y)
	$$
	that is, $X^{\pi} \in \mathcal{C}_G(u')$. Symmetrically, $\pi^{-1}$ maps $\mathcal{C}_G(u')$ into $\mathcal{C}_G(u)$. Finally, $k[X^{\pi}, Y^{\pi}; u'] = k[X,Y;u]$.
\end{proof}

We set $k_0 = \kappa(v)$ for $v \in L$ and $k_1 = \kappa(w)$ for $w \in U$.

\section{Every card have only two extensions}

In light of the results obtained in the previous section, it will be convenient to introduce new notation.
\par Let $v \in V$ and let $K \in \mathcal{C}_G(v) \setminus \{G\}$. Since the graph $G$ is fixed, let us denote $z_K := z[G,K;v]$ and we denote by $S(K) = \{x : z_K(x) = 1\}$ the set of neighbours that $v$ loses in $K$ i.e $(v,x)$ is edge in $G$ but not in $K$. It's obiviusly that only high vertices can be lost, and only low vertices can be gained. Indeed, if $v$ loses the neighbour $x$, then the degree
of $x$ drops by one; since all degrees in $K$ are again $d$ or $d+1$,
$x$ had degree $d+1$. Similarly, a gained vertex had degree $d$. Thus,
$$S(K) \subseteq N(v) \cap U, \qquad
\{x : z_K(x) = -1\} \subseteq L \setminus \bigl(N(v) \cup \{v\}\bigr).$$
Degree of $v$ is a same in $G$ and in $K$ so it loses as many neighbours as it gains. By Corollary~\ref{c3.5} this number does not
depend on $K$, and it depends on $v$ only through its degree: it equals
$k_0 \ge 1$ if $v$ is low and $k_1 \ge 1$ if $v$ is high. Thus
$$
|S(K)| = |\{x : z_K(x) = -1\}| =
\begin{cases} k_0, & v \in L,\\ k_1, & v \in U.\end{cases}
$$
For $K_1 = G$, $K_2 = K$, $s = \chi_G(v)$, Lemma~\ref{l3.1}(c) and Corollary~\ref{c3.5}):
\begin{equation}\label{eq1}\tag{$\ast$}
	(A_G - \lambda_v I)\, z_K = \kappa(v)\, \varepsilon_v,
	\qquad
	\lambda_v =
	\begin{cases}
		\theta, & v \in L,\\
		\theta-1, & v \in U.
	\end{cases}
\end{equation}
If $v \in L$ and $K_1, K_2 \in \mathcal{C}_G(v) \setminus \{G\}$, then subtracting the
two equations~\eqref{eq1} (they have the same right-hand side
$k_0\varepsilon_v$) gives

$$A_G f = \theta f \qquad \text{for } f = z_{K_1} - z_{K_2}.$$
Now, for any high vertex $w$ we fixed any some kind of extension $A_w \in \mathcal{C}_G(v) \setminus \{G\}$. We write $z_w = z_{A_w}$ and
$$
R_1(w) = \{x : z_w(x) = 1\} \quad(\text{the neighbours lost by } w),
$$
$$
S_1(w) = \{x : z_w(x) = -1\} \quad(\text{the neighbours gained by } w).
$$
Based on the above
$$
R_1(w) \subseteq N(w) \cap U, \qquad
S_1(w) \subseteq L \setminus N(w), \qquad
|R_1(w)| = |S_1(w)| = k_1 .
$$
\begin{lem}\label{l4.1}
	\medspace
	\begin{enumerate}
		\item[(a)] Let $v$ be low, $K \in \mathcal{C}_G(v)\setminus\{G\}$, and let $w$ be
		high, $K' \in \mathcal{C}_G(w)\setminus\{G\}$. Then
		$$
		\langle z_K, z_{K'} \rangle = k_1\, z_K(w) - k_0\, z_{K'}(v).
		$$
		\item[(b)] Let $v, v'$ be vertices of the same degree,
		$K \in \mathcal{C}_G(v)\setminus\{G\}$ and
		$K' \in \mathcal{C}_G(v')\setminus\{G\}$. Then
		$z_K(v') = z_{K'}(v)$. In particular, for high vertices 
		$$
		v \in R_1(v') \iff v' \in R_1(v).
		$$
	\end{enumerate}
\end{lem}

\begin{proof}
(a) By~\eqref{eq1} for $w \in U$ we have
$A_G z_{K'} = (\theta-1) z_{K'} + k_1 \varepsilon_w$, for $v \in L$ we have $A_G z_K = \theta z_K + k_0
\varepsilon_v$. Multiply on left the first equation by $z_K^{\mathsf{T}}$, and second equation by $z_K'^{\mathsf{T}}$. Since $z_K^{\mathsf{T}}z_K' = \langle z_K, z_K' \rangle$ and $A_G^{\mathsf{T}}=A_G$ we got that the two left-hand sides are equal. Equating the
right-hand sides and solving for $\langle z_K, z_{K'}\rangle$ gives
$$
(\theta-1)\langle z_K, z_{K'}\rangle + k_1\, z_K(w)=\theta \langle z_K, z_{K'}\rangle + k_0\, z_{K'}(v)
$$
hence,
$$\langle z_K, z_{K'} \rangle = k_1\, z_K(w) - k_0\, z_{K'}(v).$$
(b) Now $v$ and $v'$ have the same degree, so~\eqref{eq1} involves the
same $\lambda$ and the same $k_i \ge 1$ for both of them. The same
computation gives
$$
\lambda\langle z_K, z_{K'}\rangle + k_i\, z_K(v')
= z_K^{\mathsf T} A_G z_{K'}
= z_{K'}^{\mathsf T} A_G z_K
= \lambda\langle z_K, z_{K'}\rangle + k_i\, z_{K'}(v),
$$
so $k_i z_K(v') = k_i z_{K'}(v)$, and we may divide by $k_i$. For
high $v, v'$ have
$v \in R_1(v') \iff z_{v'}(v) = 1 \iff z_v(v') = 1 \iff v' \in R_1(v)$.
\end{proof}

\begin{lem}\label{l4.2}
	Let $A_G f = \theta f$. Then:
	\begin{enumerate}
		\item[(a)] $f(x) = 0$ for every low vertex $x$;
		\item[(b)] for every low vertex $v$ and every $K \in \mathcal{C}_G(v)\setminus\{G\}$,
		the function $f$ is constant on $S(K)$;
		\item[(c)] for every high vertex $w$, we have $f(x) = f(w)$ for all
		$x \in R_1(w)$.
	\end{enumerate}
\end{lem}

\begin{proof}
	(a) This is Lemma~\ref{l3.3} applied to $K = G$ and the class $s = 0$ of
	low vertices: an eigenvector of $A_G$ with eigenvalue $\theta - 0$ vanishes
	on all low vertices of~$G$.
	
	(b) Write $z = z_K$ and $k = k_0 \ge 1$. The matrices
	$A_K$ and $A_G$ differ only in the row and the column of $v$, we have $(A_K f)(x) = (A_G f)(x) = \theta f(x)$ if $x \neq v$. At the vertex
	$v$ itself, the row of $A_K$ is $a_K(v) = a_G(v) - z$, so
	$$
	(A_K f)(v) = \langle a_G(v) - z, f\rangle = (Af)(v) - \langle z, f\rangle
	= \theta f(v) - C = -C,
	\qquad \text{where } C := \langle z, f\rangle .
	$$
	Now apply Lemma~\ref{l3.1} to the pair $(K,G)$ (in this order) at the vertex $v$. The vertex $v$ is low in $K$
	as well (it keeps its degree), the corresponding vector is \\
	$z' = z[K,G;v] = -z$, and $k[K,G;v] = k[G,K;v] = k$. Lemma~\ref{l3.1}(c) gives
	$$
	(A_K - \theta I)\, z' = k\, \varepsilon_v .
	$$
	Put $g = f + \frac{C}{k} z'$, then
	$$
	(A_K - \theta I)\, g = -C\varepsilon_v + \tfrac{C}{k}\, k\, \varepsilon_v = 0,
	$$
	By Lemma~\ref{l3.3} applied to $K$, the function $g$ vanishes at every vertex that is low in $K$.
	\par Let $x \in S(K)$, $x$ is high in $G$
	and loses the neighbour $v$, so $x$ is low in $K$; hence $g(x) = 0$. Since
	$z'(x) = -z(x) = -1$, this reads $f(x) - \frac{C}{k} = 0$. Thus $f \equiv \frac{C}{k}$ on $S(K)$.
	
	(c) By \ref{eq1} and Lemma \ref{l3.1}(a) we have 
	$$
	 k_1 f(w)
	= f^{\mathsf T} \bigl(A_G - (\theta-1)I\bigr) z_w
	= \bigl((A_G - (\theta-1)I) f\bigr)^{\mathsf T} z_w
	= \langle f, z_w\rangle = \sum_{x \in R_1(w)} f(x)
	$$
	because $S_1(w) \subseteq L$.
	Summing by $w \in U$ and use Lemma \ref{l4.1}(b) for $R_1$, $|R_1(x)|=k_1$ we have
	$$
	\sum_{w \in U} \sum_{x \in R_1(w)} (f(x) - f(w))^2 = \sum_{x \in U} k_1 f(x)^2 - 2 \sum_{w \in U} k_1 f(w)^2 + \sum_{w \in U} k_1 f(w)^2 = 0.
	$$
	Since all terms is non zero hence they are equal to zero.
\end{proof}

\begin{lem}\label{l4.3}
	Let $v$ be low and $K_1, K_2 \in \mathcal{C}_G(v)\setminus\{G\}$. If
	$S(K_1) \cap S(K_2) \ne \varnothing$, then $K_1 = K_2$.
\end{lem}

\begin{proof}
Let $f = z_{K_1} - z_{K_2}$, then $A_Gf = \theta f$ and by Lemma \ref{l4.2}(a) $f=0$ on $L$. On $U$ vector $z_{K_i}$ is the indicator function of
$S(K_i)$, so $f = \mathbf 1_{S(K_1)} - \mathbf 1_{S(K_2)}$ there.
\par Take $x_0 \in S(K_1) \cap S(K_2)$; then $f(x_0) = 1 - 1 = 0$. By
Lemma~\ref{l4.2}(b), $f$ is constant on
$S(K_1)$, hence $f \equiv 0$ on $S(K_1)$. For $x \in S(K_1)$ this means
$\mathbf 1_{S(K_2)}(x) = \mathbf 1_{S(K_1)}(x) = 1$, i.e.\
$S(K_1) \subseteq S(K_2)$. Exchanging the roles of $K_1$ and $K_2$ gives
$S(K_2) \subseteq S(K_1)$, so $S(K_1) = S(K_2)$.

Consequently $f = 0$ on $U$ as well, so $f = 0$, i.e.\ $z_{K_1} = z_{K_2}$,
and $K_1 = K_2$
\end{proof}

\begin{thm}\label{thm4.4}
	For every vertex $v$ we have $|\mathcal{C}_G(v)| = 2$.
\end{thm}

\begin{proof}
	Let $v$ is low. Suppose, to the contrary, that there
	are two different completions $K_1 \ne K_2$ in
	$\mathcal{C}(v)\setminus\{G\}$. As in Lemma~\ref{l4.3}, the vector
	$f = z_{K_1} - z_{K_2}$ satisfies $A_Gf = \theta f$, vanishes on $L$, and equals
	$\mathbf 1_{S(K_1)} - \mathbf 1_{S(K_2)}$ on~$U$. Moreover,
	$S(K_1) \cap S(K_2) = \varnothing$ by Lemma~\ref{l4.3}.\\
	
	\textbf{Step 1:} Since $|S(K_1)| = k_0 \ge 1$, we may
	choose a vertex $w \in S(K_1)$; it is high in $G$. As $S(K_1)$ and
	$S(K_2)$ are disjoint, $w \notin S(K_2)$, so
	$$
	f(w) = \mathbf 1_{S(K_1)}(w) - \mathbf 1_{S(K_2)}(w) = 1 - 0 = 1 .
	$$
	Now let $x \in R_1(w)$. Then $x$ is high, and by Lemma~\ref{l4.2}(c)
	$f(x) = f(w) = 1$. Since $f(x) = \mathbf 1_{S(K_1)}(x) - \mathbf 1_{S(K_2)}(x)$
	is a difference of two numbers from $\{0,1\}$, the value $1$ forces
	$x \in S(K_1)$. Thus $R_1(w) \subseteq S(K_1)$. Moreover $w \notin R_1(w)$,
	because $z_w(w) = 0$. Hence
	$$
	R_1(w) \subseteq S(K_1) \setminus \{w\},
	\qquad\text{so}\qquad
	k_1 = |R_1(w)| \le |S(K_1)| - 1 = k_0 - 1 .
	$$
	
	\textbf{Step 2:} We keep the high vertex $w$ from
	Step~1 and its fixed extension $A_w$. Since
	$|S_1(w)| = k_1 \ge 1$, we may choose low $v' \in S_1(w)$ and $v'$ is not adjacent to $w$ in~$G$. The vertices $v$ and $v'$ are both
	low, so by Corollary~\ref{c3.5}(b) there is a permutation $\pi$ of $V$ with
	$\pi(v) = v'$ that maps $\mathcal{C}_G(v)$ bijectively onto
	$\mathcal{C}_G(v')$. Hence $G^\pi$, $K_1^\pi$, $K_2^\pi$ are three
	different extensions of the card $G - v'$. At most one of these three graphs equals $G$, so at least two of
	them, say $M_1 \ne M_2$, are different from $G$. Thus
	$M_1, M_2 \in \mathcal{C}_G(v')\setminus\{G\}$, and everything proved above for
	the low vertex $v$ applies equally to the low vertex $v'$. In particular,
	$|S(M_1)| = k_0$ and $S(M_1) \cap S(M_2) = \varnothing$.
	
	Note also that $z_{M_1}(w) = 0$, the vertex $v'$ can lose only
	neighbours (and $w$ is not a neighbour of $v'$) and can gain only low vertices
	(and $w$ is high). 
	\par Now we prove that $R_1(w) \cap S(M_1) = \varnothing$. $x \in R_1(w) \cap S(M_1)$. By Lemma~\ref{l4.3}, $x \notin S(M_2)$. The
	vector $f' = z_{M_1} - z_{M_2}$ satisfies $A_Gf' = \theta f'$ and equals $\mathbf 1_{S(M_1)} - \mathbf 1_{S(M_2)}$ on $U$. Since $x$ is high,
	$f'(x) = 1 - 0 = 1$. As $x \in R_1(w)$, Lemma~\ref{l4.2}(c) gives
	$f'(w) = f'(x) = 1$. On the other hand, $z_{M_1}(w) = 0$, so $f'(w) = -z_{M_2}(w) \in \{0,-1\}$ is contradiction.
	\par Apply Lemma~\ref{l4.1}(a) to the low vertex $v'$ with the extension $M_1$ and
	to the high vertex $w$ with the extension $A_w$. Since $z_{M_1}(w) = 0$
	and $z_w(v') = -1$ (because $v' \in S_1(w)$), we get
	$$
	\langle z_{M_1}, z_w\rangle = k_1 \cdot 0 - k_0 \cdot (-1) = k_0 .
	$$
	Now we estimate the same sum $\sum_x z_{M_1}(x)\, z_w(x)$ term by term.
	\begin{itemize}
		\item \textit{High $x$.} Both factors lie in $\{0,1\}$, and the product is
		nonzero only if $x \in S(M_1) \cap R_1(w)$. This set is empty, so all these terms vanish.
		\item \textit{Low $x$.} Both factors lie in $\{-1,0\}$, so the product is $0$
		or $1$, and it can equal $1$ only if $z_w(x) = -1$, i.e.\
		$x \in S_1(w)$. For $x = v' \in S_1(w)$ the product is $0$, because $z_{M_1}(v_2) = 0$. Hence the low terms contribute at most
		$|S_1(w) \setminus \{v'\}| = k_1 - 1$.
	\end{itemize}
	Therefore $k_0 = \langle z_{M_1}, z_w\rangle \le k_1 - 1$.
	\par Steps 1 and 2 contradict each other, so $|\mathcal{C}_G(v)| = 2$ for every low
	vertex~$v$.
	
	For $v \in U$, we apply the proven result to the counterexample $(\overline{G}, \overline{H})$ (Lemma~\ref{l2.5}), in which $v$ is low: its extensions of the card $\overline{G} - u$ are precisely the complements of the graphs in $\mathcal{C}_G(v)$.
\end{proof}

\begin{cor}\label{c4.5}
	$ \mathcal{C}_G(u) = \{G, \mathrm{alt}(u)\} $ and $ \mathrm{alt}(u) \cong H $. For vertices $ u, u' $ of the same degree, there exists $ \pi \in \operatorname{Aut}(G) $ with $ \pi(u) = u' $; thus, $ \operatorname{Aut}(G) $ is transitive on $ L $ and on $ U $. For any $\pi \in \operatorname{Aut}(G)$, we have $\mathrm{alt}(\pi u) = \mathrm{alt}(u)^\pi$, and for $z_u := z[G, \mathrm{alt}(u); u]$ it holds that $z_{\pi u}(\pi x) = z_u(x)$.
\end{cor}

\begin{proof}
	The set $\mathcal{C}(u)$ contains a graph isomorphic to $H$ (Lemma \ref{l3.3}); it is distinct from $G$, hence it is $\mathrm{alt}(u)$. For $\pi$ from Corollary \ref{c3.5}(b), the graph $G^\pi \in \mathcal{C}(u')$ is isomorphic to $G$, while $\mathrm{alt}(u') \cong H \not\cong G$; therefore $G^\pi = G$, i.e., $\pi \in \operatorname{Aut}(G)$. If $\pi \in \operatorname{Aut}(G)$, then, as in Corollary \ref{c3.5}(b), $\mathrm{alt}(u)^\pi \in \mathcal{C}(\pi u)$ and $\mathrm{alt}(u)^\pi \neq G^\pi = G$, so that $\mathrm{alt}(u)^\pi = \mathrm{alt}(\pi u)$. The last equality follows by substitution.
\end{proof}

\section{Main identities}

Now, due to the fact that we have only two extensions, we will fix $z_u = z[G,\alt(u);u]$ for all $u \in V$ and

\begin{equation}\label{eq2}\tag{$\ast \ast$}
	(A_G-\lambda_u I)\,z_u = k_u\,\varepsilon_u,
	\qquad
	 k_u=\begin{cases} k_0, & v \in L,\\ k_1, & v \in U.\end{cases}
\end{equation}
and $\lambda_u$ same as in previous section.
\par For $v \in L$, $w \in U$ denote:
\begin{itemize}
	\item $S_0(v)=\{x : z_v(x)=1\}\subseteq N(v) \cap U$
	\item $R_0(v)=\{x : z_v(x)=-1\} \subseteq L \setminus (N(v) \cup \{v\})$
	\item $R_1(w)=\{x : z_w(x)=1\} \subseteq N(w) \cap U$
	\item $S_1(w)=\{x : z_w(x)=-1\} \subseteq L \setminus N(w)$
\end{itemize}

We now introduce $0/1$-matrices. The rows and columns of $R_0$, $E_0$ are
indexed by the vertices of $L$; those of $R_1$, $E_1$ by the vertices of $U$;
the matrices $D_0,D_1,D_2,D_3$ have rows indexed by $v\in L$ and columns indexed
by $w\in U$:
$$
\begin{aligned}
	R_0(v,v') &= [{v'\in R_0(v)}], &\qquad R_1(w,w') &= [{w'\in R_1(w)}],\\
	D_0(v,w) &= [{w\in S_0(v)}], &\qquad D_1(v,w) &= [{v\in S_1(w)}],\\
	D_2(v,w) &= [{v\sim w,\ w\notin S_0(v)}], &\qquad D_3(v,w) &= [{v\nsim w,\ v\notin S_1(w)}],\\
	E_0(v,v') &= [{v'\neq v,\ v'\nsim v,\ v'\notin R_0(v)}], &\qquad E_1(w,w') &= [{w'\sim w,\ w'\notin R_1(w)}].
\end{aligned}
$$

\textbf{Remark:} $D_0 + D_1 + D_2 + D_3 = J$ where $J$ is the all-ones matrix of
size $p\times q$.
\begin{lem}\label{l5.1}
	\medspace
	\begin{enumerate}
		\item[(a)] The matrices $R_0$ and $R_1$ are symmetric:
		$v'\in R_0(v)\iff v\in R_0(v')$ and $w'\in R_1(w)\iff w\in R_1(w')$.
		\item[(b)] For $v\in L$, $w\in U$ and all $x$ one has $z_v(x)\,z_w(x)\ge 0$, and
		$$
			|S_0(v)\cap R_1(w)|+|R_0(v)\cap S_1(w)| = \langle z_v, z_w\rangle
			= k_1[{w\in S_0(v)}]+k_0[{v\in S_1(w)}].
		$$
		In matrix form,
		
			$$R_0D_1+D_0R_1 = k_1D_0+k_0D_1.$$
		
		\item[(c)] Every row of $D_0$ contains exactly $k_0\ge 1$ ones, and every column of
		$D_1$ contains exactly $k_1$ ones. All column sums of $D_0$ are equal to
		one and the same number~$\gamma$, and all row sums of $D_1$ are equal to
		one and the same number~$\delta$; moreover, $pk_0=q\gamma$,
		$qk_1=p\delta$, $\gamma\delta=k_0k_1$ and $\delta\ge 1$. In particular,
		every low vertex lies in some set $S_1(w)$.
	\end{enumerate}
\end{lem}

\begin{proof}
	(a) By Lemma \ref{l4.1}(b) $z_{u'}(u)=z_u(u')$ for all vertices $u,u'$ of the same degree. If $u,u' \in L$ it's proves symetry for $R_0$ and similarly for $u,u' \in U$.
	\par(b) On $L$ both vectors $z_v,z_w$ take values in $\{-1,0\}$, and on $U$ in
	$\{0,1\}$. Hence at every point $x$ their values do not have opposite signs, and
	$z_v(x)z_w(x)\in\{0,1\}$. The product equals $1$ in exactly two cases: $x\in U$
	and $z_v(x)=z_w(x)=1$, i.e.\ $x\in S_0(v)\cap R_1(w)$; or $x\in L$ and
	$z_v(x)=z_w(x)=-1$, i.e.\ $x\in R_0(v)\cap S_1(w)$. Summing over $x$, we obtain
	$$
	\langle z_v, z_w\rangle=|S_0(v)\cap R_1(w)|+|R_0(v)\cap S_1(w)|.
	$$
	By Lemma \ref{l4.1}(a) we have $\langle z_v, z_w\rangle = k_1z_v(w)-k_0z_w(v)$. Since
	$w\in U$, then $z_v(w)=\mathbf 1_{S_0(v)}(w)=\mathbf 1{w\in S_0(v)}$. Since
	$v\in L$, then $z_w(v)=-\mathbf 1_{S_1(w)}(v)=-\mathbf 1{v\in S_1(w)}$.
	Substituting, we obtain the second equality.\\
	For matrix form we compute the $(v,w)$ entry of each term on the left-hand side:
	$$
	\begin{aligned}
		(R_0D_1)(v,w) &= \sum_{v'\in L}[{v'\in R_0(v)}]\,[{v'\in S_1(w)}]
		= |R_0(v)\cap S_1(w)|,\\
		(D_0R_1)(v,w) &= \sum_{w'\in U}[{w'\in S_0(v)}]\,[{w\in R_1(w')}]
		= \sum_{w'\in U}[{w'\in S_0(v)}]\,[{w'\in R_1(w)}]
		= |S_0(v)\cap R_1(w)|;
	\end{aligned}
	$$
	The $(v,w)$ entry of the right-hand side equals
	$k_1[{w\in S_0(v)}]+k_0[{v\in S_1(w)}]$. Hence matrix form is precisely
	identity written simultaneously for all pairs
	$(v,w)\in L\times U$.
	\par(c) \textit{Rows of $D_0$ and columns of $D_1$.}
	The sum of the row $v$ of $D_0$ equals
	$\#\{w\in U : w\in S_0(v)\}=|S_0(v)|=k_0$ (recall that $S_0(v)\subseteq U$), and
	$k_0\ge 1$ by Lemma~\ref{l3.1}(d). The sum of the column $w$ of $D_1$ equals
	$\#\{v\in L : v\in S_1(w)\}=|S_1(w)|=k_1$.
	
	\textit{Columns of $D_0$ and rows of $D_1$.}
	Let $\pi\in\Aut(G)$. Automorphisms preserve degrees, so $\pi(L)=L$ and
	$\pi(U)=U$. By Corollary~\ref{c4.5}, $z_{\pi(v)}(\pi (x))=z_v(x)$ for all $v,x$; hence
	$$
	w\in S_0(v)\iff \pi (w)\in S_0(\pi (v)),\qquad
	v\in S_1(w)\iff \pi (v)\in S_1(\pi (w)).
	$$
	Take $w,w'\in U$. By Corollary~\ref{c4.5} there exists $\pi\in\Aut(G)$ with
	$\pi (w)=w'$. The map $v\mapsto\pi (v)$ is a bijection $L\to L$, and by the first
	equivalence it maps the set $\{v\in L : w\in S_0(v)\}$ onto the set
	$\{v''\in L : w'\in S_0(v'')\}$. Therefore the columns $w$ and $w'$ of $D_0$
	contain the same number of ones; we denote this number by~$\gamma$. In the same
	way, for $v,v'\in L$ and $\pi\in\Aut(G)$ with $\pi (v)=v'$, the map
	$w\mapsto\pi (w)$ maps $\{w\in U : v\in S_1(w)\}$ onto
	$\{w''\in U : v'\in S_1(w'')\}$; hence all row sums of $D_1$ are equal to one
	and the same number~$\delta$.
	\par The number of ones in $D_0$ equals $pk_0$ (summing over the $p$ rows) and also
	$q\gamma$ (summing over the $q$ columns), so $pk_0=q\gamma$. The number of ones
	in $D_1$ equals $qk_1$ (by columns) and $p\delta$ (by rows), so $qk_1=p\delta$.
	Multiplying these identities yields $pq\,k_0k_1=pq\,\gamma\delta$, and since
	$p,q\ge 1$ we get $\gamma\delta=k_0k_1$. Finally, $\delta=\frac{qk_1}{p}>0$ and $\delta$
	is an integer, so $\delta\ge 1$. This means that every row $v$ of $D_1$
	contains a one, i.e.\ $v\in S_1(w)$ for some $w\in U$.
\end{proof}

\begin{lem}\label{l5.2}
	The following identity holds:
	$$
	k_1-k_0=\theta.
	$$
\end{lem}
 \begin{proof}
 	\medspace
 	\begin{enumerate}
 		\item Take any $w\in U$ (possible since $q\ge 1$) and any $v\in S_1(w)$ (possible
 		since $|S_1(w)|=k_1\ge 1$). Then $v\in L$ and $v\nsim w$ and $S_0(v)\subseteq N(v)$, we have $w\notin S_0(v)$; moreover,
 		$v\in S_1(w)$. By Lemma \ref{l5.1}(b),
 		$$
 		\langle z_v, z_w\rangle=k_1\cdot 0+k_0\cdot 1=k_0 .
 		$$
 		By \ref{eq2}
 		$$
 		(A_G-\lambda_{w}I)z_w(v)=k_w\varepsilon_w(v) \qquad z_w(v)=1
 		$$
 		Hence $(A_Gz_w)(v)=(\theta-1)z_w(v)=1-\theta$.
 		\par Put $a'=a_{\alt(v)}(v)=a_G(v)-z_v$. Since $(A_Gz_w)(v)=\langle a_G(v), z_w\rangle$,
 		$$
 		\sum_y a'(y)\,z_w(y)=\langle a_G(v), z_w\rangle-\langle z_v, z_w\rangle=(1-\theta)-k_0 .
 		$$
 		We have $a'(y)\in\{0,1\}$ and $z_w(y)\in\{-1,0,1\}$. If $y\notin S_1(w)$, then
 		$z_w(y)\ge 0$ and the summand $a'(y)z_w(y)$ is nonnegative. If $y\in S_1(w)$,
 		then the summand is at least $-1$. Moreover, $v\in S_1(w)$ but $a'(v)=0$, so the summand with $y=v$ vanishes. Therefore at
 		most $|S_1(w)|-1=k_1-1$ summands can be negative, each of them being at least
 		$-1$:
 		$$
 		\sum_y a'(y)\,z_w(y)\ge -(k_1-1).
 		$$
 		Comparing with the exact value, we get $1-\theta-k_0\ge 1-k_1$, that is,
 		$k_1\ge k_0+\theta$.
 		\item Take any $v\in L$ and any $w\in S_0(v)$.
 		Then $w\in U$ and $v\sim w$.
 		Since $v\sim w$ and $S_1(w)\cap N(w)=\varnothing$, we have $v\notin S_1(w)$. By Lemma \ref{l5.1}(b)
 		$$
 		\langle z_v, z_w\rangle=k_1\cdot 1+k_0\cdot 0=k_1 .
 		$$
 		By \ref{eq2} $(A_Gz_v)(w)=\theta z_v(w)=\theta$
 		\par Put $a''=a_{\alt(w)}(w)=a_G(w)-z_w$. Then
 		$$
 		\sum_y a''(y)\,z_v(y)=\langle a_G(w), z_w\rangle-\langle z_v, z_w\rangle=\theta-k_1.
 		$$
 		Since $a''(y)\in\{0,1\}$, the summand $a''(y)z_v(y)$ is nonnegative if
 		$z_v(y)\ge 0$, and is at least $-1$ if $z_v(y)=-1$, i.e.\ if $y\in R_0(v)$.
 		Hence every summand is at least $-[{y\in R_0(v)}]$, and
 		$$
 		\sum_y a''(y)\,z_v(y)\ge -|R_0(v)|=-k_0 .
 		$$
 		Consequently $\theta-k_1\ge -k_0$, that is, $k_1\le k_0+\theta$.
 		
 		\medskip
 		Steps~1 and~2 together give $k_1=k_0+\theta$
 	\end{enumerate}
 \end{proof}
 
 \begin{prop}\label{p5.3}
 	The following identity holds:
 	\begin{enumerate}
 		\item[(a)] $R_0D_3 + D_0E_1 = k_0(D_2 + D_3)$
 		\item[(b)] $E_0D_1 + D_2R_1 = k_1(D_2 + D_3)$
 	\end{enumerate}
 \end{prop}
 \begin{proof}
 	Fix $v \in L$ and $w \in U$, and set
 	$\sigma_0 = [\, w \in S_0(v) \,]$ and $\sigma_1 = [\, v \in S_1(w) \,]$;
 	then $(D_2 + D_3)(v,w) = 1 - \sigma_0 - \sigma_1$.
 	
 	(a)We have $(R_0D_3)(v,w) = \bigl|R_0(v) \setminus \bigl(N(w) \cup S_1(w)\bigr)\bigr|$
 	and, by the symmetry of $R_1$,
 	$(D_0E_1)(v,w) = \bigl|\bigl(S_0(v) \cap N(w)\bigr) \setminus R_1(w)\bigr|$.
 	Since $S_1(w) \cap N(w) = \varnothing$ and $R_1(w) \subseteq N(w)$, it follows that
 	\[
 	(R_0D_3 + D_0E_1)(v,w)
 	= k_0 - |R_0(v) \cap N(w)| - |R_0(v) \cap S_1(w)|
 	+ |S_0(v) \cap N(w)| - |S_0(v) \cap R_1(w)|.
 	\]
 	Here $|S_0(v) \cap N(w)| - |R_0(v) \cap N(w)| = (A_Gz_v)(w) = \theta\sigma_0$ by~\ref{eq2},
 	while the sum of the remaining two intersections equals $k_1\sigma_0 + k_0\sigma_1$ by Lemma \ref{l5.1}(b).
 	By Lemma \ref{l5.2} we obtain
 	$k_0 + \theta\sigma_0 - k_1\sigma_0 - k_0\sigma_1 = k_0(1 - \sigma_0 - \sigma_1)$.
 	
 	(b)We have $(E_0D_1)(v,w) = \bigl|S_1(w) \setminus \bigl(N(v) \cup \{v\} \cup R_0(v)\bigr)\bigr|$
 	and $(D_2R_1)(v,w) = \bigl|\bigl(R_1(w) \cap N(v)\bigr) \setminus S_0(v)\bigr|$.
 	Since the sets $N(v)$, $\{v\}$, $R_0(v)$ are pairwise disjoint and $S_0(v) \subseteq N(v)$, we get
 	\[
 	(E_0D_1 + D_2R_1)(v,w)
 	= k_1 - |S_1(w) \cap N(v)| - \sigma_1 - |S_1(w) \cap R_0(v)|
 	+ |R_1(w) \cap N(v)| - |R_1(w) \cap S_0(v)|.
 	\]
 	By \ref{eq2},
 	$|R_1(w) \cap N(v)| - |S_1(w) \cap N(v)| = (A_Gz_w)(v) = (\theta - 1)z_w(v) = -(\theta - 1)\sigma_1$,
 	and by Lemma \ref{l5.1}(b), by Lemma \ref{l5.2}, we obtain
 	$k_1 - (\theta - 1)\sigma_1 - \sigma_1 - k_1\sigma_0 - k_0\sigma_1 = k_1(1 - \sigma_0 - \sigma_1)$.
 \end{proof}
 
 \section{Connectivity of switchig graph}
 In this section, we discuss some general theoretical properties; we do not use the notation or the assumption regarding a counterexample unless stated otherwise.
\par A graph \textit{module} of a graph $G=(V;E)$ is a vertex set $B \subseteq V$ such that every vertex outside $B$ is adjacent either to every vertex in $B$ or to no vertex in $B$. Equivalently, vertices in $B$ have identical neighbors outside $B$. 
\par For non empity module $B$ and $y \in V \setminus B$ denote $\beta_G(y)=1$ in first case and $\beta_G(y)=0$ another. Thus $G(b,y)=\beta_G(y)$ for all $b \in B, y \notin B$. If $B \neq V$ than $B$ is called \textit{proper} module.
\par Let $F=(W;E'), G=(V;E)$ - graphs. Injective map $f: W \to V$ is called \textit{embedding} $F$ to $G$, if for all $a,a' \in W$
$$a \sim_F a' \iff f(a) \sim_G f(a')$$
Put $\emb(F,G)$ is number of all such embeddings and $\emb_X(F,G)$ for $X \subseteq V$ is number of such embeddings with condition $f(W) \subseteq X$. For $S \subseteq W$ we denoted induced subgraph by $F[S]$.
\par Not hard to see that $\emb_{V\setminus\{v\}}(F,G)=\emb(F,G - v)$ and if $G \cong G'$ than $\emb(F,G)=\emb(F,G')$.
\par Denote by $s(F,G)$ the number of copies of $F$ in $G$. It's easy to see that $\emb(F,G)=|\Aut(F)|s(F,G)$.
\begin{lem}\label{l6.1}
	Let $G\sim K$ on a set $V$, where $|V| = n$ and let $F$ be a graph on a set $W$ with $|W|<n$. Then $\emb(F,G)=\emb(F,K)$.
\end{lem}
\begin{proof}
	By $\emb(F,G)=|\Aut(F)|s(F,G)$ and \cite[Lemma 2.20]{BonMur2008} we have 
	$$(n-|W|)\,\emb(F,G)=\sum_{v\in V}\emb(F,G-v).$$
	Since $G \sim K$ we have a bijection $\sigma$, such that $G-v \cong K - \sigma(v)$. Hence $$\emb(F,G-v)=\emb(F,K-\sigma(v))$$ and 
	$$(n-|W|)\,\emb(F,G)=(n-|W|)\,\emb(F,K)$$
	since $(n-|W|)\ge 0$ we have $\emb(F,G)=\emb(F,K).$
\end{proof}
For $S\subseteq W$ denote by $\out(S;G)$ the number of injective maps $g\colon W\setminus S\to V\setminus B$ which are induced embeddings of $F[W\setminus S]$ into $G$
\begin{lem}[Decomposition]\label{l6.2}
	Let $G$ be a graph on $V$, $B$ a nonempty module of $G$, and $F$ a graph on a vertex set $W$ and satisfies condition
	$$ \beta_G(g(c))=F(a,c)\qquad\text{for all}\ a\in S,\ c\in W\setminus S .$$
	Then
	$$ \emb(F,G)=\sum_{S\subseteq W}\emb_B(F[S],G)\cdot\out(S;G).$$
\end{lem}
\begin{proof}
	To every induced embedding $f$ of $F$ into $G$ we assign the triple $(S,f_1,f_2)$, where $S=f^{-1}(B)$, $f_1=f|_S$ and $f_2=f|_{W\setminus S}$. Then $f_1$ is an induced embedding of $F[S]$ into $G$ with image in $B$, and $f_2$ is an induced embedding of $F[W\setminus S]$ into $M$ with image in $V\setminus B$. Moreover, for $a\in S$, $c\in W\setminus S$ we have $f(a)\in B$, $f(c)\notin B$, and since $B$ is a module, $G(f(a),f(c))=\beta_G(f(c))$; since $f$ is an induced embedding, this number equals $F(a,c)$. Hence $f_2$ satisfies lemma condition.
	Conversely, suppose we are given $S\subseteq W$, an induced embedding $f_1$ of $F[S]$ with image in $B$, and a map $f_2$ as in the definition of $\out(S;M)$. The glued map $f\colon W\to V$ ($f=f_1$ on $S$, $f=f_2$ on $W\setminus S$) is injective, since $f_1$ and $f_2$ are injective and their images lie in the disjoint sets $B$ and $V\setminus B$. It is an induced embedding: for pairs inside $S$ and inside $W\setminus S$ this is the corresponding property of $f_1$ and $f_2$, and for $a\in S$, $c\in W\setminus S$ we have $G(f(a),f(c))=\beta_G(f_2(c))=F(a,c)$ by the module property and lemma condition. Moreover, $f^{-1}(B)=S$.
	The two constructions are mutually inverse. Therefore, for fixed $S$, the number of embeddings $f$ with $f^{-1}(B)=S$ equals $\emb_B(F[S],G)\cdot\out(S;G)$, and summing over $S$ gives a necessary equality.
\end{proof}
We emphasize that the number $\out(S;G)$ is determined by the induced subgraph $G[V\setminus B]$ and by the function $\beta_G$ on $V\setminus B$, that is, only by the adjacencies of $G$ on pairs of vertices not contained in $B$.

	\begin{lem}\label{l6.3}
			Let $G\sim K$, let $B$ be a nonempty module of the graph $G$, and suppose that the graphs $G$ and $K$ coincide on all pairs of vertices not contained in $B$. Then
		\begin{itemize}
			\item[(a)] $B$ is a module of the graph $K$, $\beta_K=\beta_G$, and $\out(S;K)=\out(S;G)$ for every graph $F$ on $W$ and every $S\subseteq W$;
			\item[(b)] $\emb_B(F,G)=\emb_B(F,K)$ for every graph $F$ on a set $W$ with $|W|<n$.
		\end{itemize}
	\end{lem}
	\begin{proof}
		(a) For $b\in B$ and $y\notin B$ the pair $\{b,y\}$ is not contained in $B$, hence $K(b,y)=G(b,y)=\beta_G(y)$, that is, $B$ is a module of $K$ and $\beta_K=\beta_G$. Pairs of vertices of $V\setminus B$ are not contained in $B$ either, hence $K[V\setminus B]=G[V\setminus B]$. By the remark above, $\out(S;K)=\out(S;G)$.
		\par (b) Induction on $|W|$. Separating in Lemma \ref{l6.2} the term with $S=W$ (in which $\out(W;M)=1$), we obtain for $M\in\{G,K\}$
		$$ \emb_B(F,M)=\emb(F,M)-\sum_{S\subsetneq W}\emb_B(F[S],M)\,\out(S;M). $$
		Since $|W|<n$, Lemma~\ref{l6.1} gives $\emb(F,G)=\emb(F,K)$. In every term of the sum $|S|<|W|<n$, so by the induction hypothesis $\emb_B(F[S],G)=\emb_B(F[S],K)$, while $\out(S;G)=\out(S;K)$ by (a). Hence $\emb_B(F,G)=\emb_B(F,K)$. The base of the induction ($W=\varnothing$) is contained in this argument: the sum is empty and $\emb(F,M)=1$.
	\end{proof}
	\begin{thm}\label{t6.4}
		Let graphs $G$ and $K$ on $V$ have the same deck, let $B\subsetneq V$ be a module of the graph $G$, and suppose that $G$ and $K$ coincide on all pairs of vertices not contained in $B$. Then $G\cong K$.
	\end{thm}
	\begin{proof}
		If $B=\varnothing$, then $G$ and $K$ coincide on all pairs, that is, $G=K$. Let $B\ne\varnothing$. Apply Lemma~\ref{l6.3}(b) to the graph $F=G[B]$ on the set $W=B$; this is correct since $|B|<n$. The identity inclusion $B\hookrightarrow V$ is an induced embedding of $G[B]$ into $G$ with image in $B$, hence
		$$ \emb_B(G[B],K)=\emb_B(G[B],G)\ge 1. $$
		Therefore there is an injective map $\varphi\colon B\to B$ such that $x\sim_G y\iff\varphi(x)\sim_K\varphi(y)$ for all distinct $x,y\in B$.
		Extend $\varphi$ to a bijection $\Phi\colon V\to V$ by setting $\Phi(y)=y$ for $y\notin B$, and let us check that $\Phi$ is an isomorphism $G\to K$, that is, $K(\Phi x,\Phi y)=G(x,y)$ for all distinct $x,y$.
		\begin{itemize}
			\item[$\bullet$] $x,y\in B$: this is the property of $\varphi$.
			\item[$\bullet$] $x,y\notin B$: the pair $\{x,y\}$ is not contained in $B$, and $K(x,y)=G(x,y)$.
			\item[$\bullet$] $x\in B$, $y\notin B$: the pair $\{\varphi(x),y\}$ is not contained in $B$, hence $K(\varphi(x),y)=G(\varphi(x),y)$; since $B$ is a module of $G$ and $\varphi(x),x\in B$, we have $G(\varphi(x),y)=\beta_G(y)=G(x,y)$.
		\end{itemize}
		Thus $G\cong K$.
	\end{proof}
	
	From this point on, $(G,H)$ serves as a counterexample once again, and all the notation from the preceding sections applies.
	\par Denote $f_L:=f|_L$ and $f_U:=f|_U$
	\begin{lem}\label{l6.5}
		Let $f\colon V\to\mathbb{Q}$ and $\sum_{u\in V}f(u)\,z_u=0$. Then
		$$ R_0f_L=k_0f_L,\qquad R_1f_U=k_1f_U,\qquad D_1f_U=-k_0f_L,\qquad D_0^{\mathsf T}f_L=-k_1f_U. $$
	\end{lem}
	\begin{proof}
		By \ref{eq2}, $(A_G-\theta I)z_u=(\lambda_u-\theta)z_u+k_u\varepsilon_u$, where $\lambda_u-\theta=0$ for $u\in L$ and $\lambda_u-\theta=-1$ for $u\in U$. Applying $A_G-\theta I$ to the equality $\sum_u f(u)z_u=0$, we obtain
		$$ -\sum_{u\in U}f(u)\,z_u+\sum_{u\in V}k_u f(u)\,\varepsilon_u=0. $$
		We write this equality in the coordinate $y\in V$, setting $k_y=k_0$ for $y\in L$ and $k_y=k_1$ for $y\in U$:
		\[ \sum_{u\in U}f(u)\,z_u(y)=k_y f(y). \tag{a} \]
		Subtracting (a) from the equality $\sum_{u\in V}f(u)z_u(y)=0$, we get
		\[ \sum_{u\in L}f(u)\,z_u(y)=-k_y f(y). \tag{b} \]
		We compute the left-hand sides of (a) and (b) using main identities.
		\begin{itemize}
			\item (a) for $y\in U$: for $u\in U$ we have $z_u(y)=[\,y\in R_1(u)\,]=R_1(u,y)=R_1(y,u)$ by the symmetry of $R_1$, so the left-hand side equals $(R_1f_U)(y)$. We obtain $R_1f_U=k_1f_U$.
			\item (a) for $y\in L$: for $u\in U$ we have $z_u(y)=-[\,y\in S_1(u)\,]=-D_1(y,u)$, so the left-hand side equals $-(D_1f_U)(y)$. We obtain $D_1f_U=-k_0f_L$.
			\item (b) for $y\in L$: for $u\in L$ we have $z_u(y)=-[\,y\in R_0(u)\,]=-R_0(y,u)$, so the left-hand side equals $-(R_0f_L)(y)$. We obtain $R_0f_L=k_0f_L$.
			\item (b) for $y\in U$: for $u\in L$ we have $z_u(y)=[\,y\in S_0(u)\,]=D_0(u,y)$, so the left-hand side equals $(D_0^{\mathsf T}f_L)(y)$. We obtain $D_0^{\mathsf T}f_L=-k_1f_U$.
		\end{itemize}
	\end{proof}
	\begin{lem}\label{l6.6}
		\medspace
		\begin{itemize}
			\item[(a)] Let $D$ be a $0/1$ matrix indexed by $L\times U$ with exactly $r$ ones in every row and exactly $c$ ones in every column. Then for $\varphi\colon L\to\mathbb{Q}$, $h\colon U\to\mathbb{Q}$
			$$ \mathcal{E}_D(\varphi,h):=\sum_{v\in L}\sum_{w\in U}D(v,w)\bigl(h(w)-\varphi(v)\bigr)^2=c\,\|h\|^2-2\langle\varphi,Dh\rangle+r\,\|\varphi\|^2. $$
			\item[(b)] Let $R$ be a symmetric $0/1$ matrix indexed by $X\times X$ with exactly $k$ ones in every row. Then for $\psi\colon X\to\mathbb{Q}$
			$$ \sum_{x,y\in X}R(x,y)\bigl(\psi(x)-\psi(y)\bigr)^2=2\bigl(k\,\|\psi\|^2-\langle\psi,R\psi\rangle\bigr). $$
		\end{itemize}
	\end{lem}
	\begin{proof}
		(a) Expand the square: $\sum_{v,w}D(v,w)h(w)^2=\sum_w c\,h(w)^2$, $\sum_{v,w}D(v,w)\varphi(v)^2=\sum_v r\,\varphi(v)^2$, $\sum_{v,w}D(v,w)\varphi(v)h(w)=\langle\varphi,Dh\rangle$. (b) Similarly: $\sum_{x,y}R(x,y)\psi(x)^2=k\|\psi\|^2$, and by the symmetry of $R$ the sum $\sum_{x,y}R(x,y)\psi(y)^2$ also equals $k\|\psi\|^2$; the cross term equals $-2\langle\psi,R\psi\rangle$.
	\end{proof}
	\par The \textit{switching graph} $\mathcal{T}$ is the simple graph on the set $V$ in which distinct vertices $u$ and $x$ are adjacent if and only if $z_u(x)\ne 0$ or $z_x(u)\ne 0$.
	\par From this definition we immediately have that pair of distinct vertices $\{u,x\}$ is an edge of $\mathcal{T}$ if and only if it is of one of the following four types:
	\begin{itemize}
		\item[$(R_0)$] $u,x\in L$ and $x\in R_0(u)$ (equivalently, $u\in R_0(x)$);
		\item[$(R_1)$] $u,x\in U$ and $x\in R_1(u)$ (equivalently, $u\in R_1(x)$);
		\item[$(D_0)$] one of the vertices, $v$, lies in $L$, the other one, $w$, lies in $U$, and $w\in S_0(v)$;
		\item[$(D_1)$] one of the vertices, $v$, lies in $L$, the other one, $w$, lies in $U$, and $v\in S_1(w)$.
	\end{itemize}
	
	\begin{lem}\label{l6.7}
		Let $f\colon V\to\mathbb{Q}$ and $\sum_{u\in V}f(u)\,z_u=0$. Define $g\colon V\to\mathbb{Q}$ by
		\[ g(v)=f(v)\ \ (v\in L),\qquad g(w)=-\frac{\delta}{k_0}\,f(w)\ \ (w\in U). \]
		Then $g(u)=g(x)$ for every edge $\{u,x\}$ of the graph $T$. Consequently, $g$ is constant on every connected component of $T$.
	\end{lem}
	\begin{proof}
		Set $\varphi=g_L=f_L$ and $h=g_U=-(\delta/k_0)f_U$. From Lemma \ref{l6.5} and the equality $\gamma\delta=k_0k_1$ we get
		\[ R_0\varphi=k_0\varphi,\qquad R_1h=k_1h, \tag{a} \]
		\[ D_1h=-\frac{\delta}{k_0}D_1f_U=\delta\varphi,\qquad D_0^{\mathsf T}\varphi=-k_1f_U=\frac{k_0k_1}{\delta}h=\gamma h. \tag{b} \]
		By Lemma~\ref{l6.6}(a) for $D=D_0$ (rows with $k_0$ ones, columns with $\gamma$ ones) and (b),
		$$ \mathcal{E}_{D_0}(\varphi,h)=\gamma\|h\|^2-2\langle D_0^{\mathsf T}\varphi,h\rangle+k_0\|\varphi\|^2=k_0\|\varphi\|^2-\gamma\|h\|^2. $$
		For $D=D_1$ (rows with $\delta$ ones, columns with $k_1$ ones) we similarly obtain, using (b),
		$$ \mathcal{E}_{D_1}(\varphi,h)=k_1\|h\|^2-2\langle\varphi,D_1h\rangle+\delta\|\varphi\|^2=k_1\|h\|^2-\delta\|\varphi\|^2. $$
		Consequently,
		$$ \delta\,\mathcal{E}_{D_0}(\varphi,h)+k_0\,\mathcal{E}_{D_1}(\varphi,h)=(k_0k_1-\gamma\delta)\|h\|^2=0. $$
		Both $\mathcal{E}_{D_0}$ and $\mathcal{E}_{D_1}$ are sums of nonnegative terms, and $\delta,k_0>0$; hence $\mathcal{E}_{D_0}=\mathcal{E}_{D_1}=0$, and every term vanishes: $h(w)=\varphi(v)$, i.e.\ $g(w)=g(v)$, whenever $D_0(v,w)=1$ or $D_1(v,w)=1$.
		\par By Lemma~\ref{l6.6}(b) and (b),
		$$ \sum_{v,v'\in L}R_0(v,v')\bigl(\varphi(v)-\varphi(v')\bigr)^2=2\bigl(k_0\|\varphi\|^2-\langle\varphi,k_0\varphi\rangle\bigr)=0, $$
		and similarly $\sum_{w,w'\in U}R_1(w,w')\bigl(h(w)-h(w')\bigr)^2=2\bigl(k_1\|h\|^2-\langle h,k_1h\rangle\bigr)=0$. Hence $g(v)=g(v')$ whenever $R_0(v,v')=1$, and $g(w)=g(w')$ whenever $R_1(w,w')=1$.
		\par Every edge of $\mathcal{T}$ is of one of the types $R_0$, $R_1$, $D_0$, $D_1$, and by regarding what was written above the values of $g$ at its endpoints coincide. Therefore $g$ is constant along every path in $\mathcal{T}$, that is, on every connected component.
	\end{proof}
	
	\begin{thm}\label{t6.8}
		The switching graph $\mathcal T$ is connected.
	\end{thm}
	
	\begin{proof}
			Let $B$ be the vertex set of some connected component of the graph $\mathcal T$; $B\ne\varnothing$. We prove that $B=V$.
		\par \textbf{Step 1:} If $u\in B$, then $\{x: z_u(x) \neq 0\} \subseteq N_T(u)\subseteq B$, since the neighbours of a vertex lie in its connected component.
		Consider the square matrix $Z_B=\bigl(z_u(x)\bigr)_{x\in B,\,u\in B}$ of order $|B|$ over $\mathbb{Q}$ (rows are indexed by the vertices $x$, columns by the vertices $u$).
		\par \textbf{Step 2:} the right kernel of $Z_B$ is at most one-dimensional. Let $c\colon B\to\mathbb{Q}$ with $Z_Bc=0$, i.e.\ $\sum_{u\in B}c(u)z_u(x)=0$ for all $x\in B$. Extend $c$ by zero to a function $f$ on $V$. The vector $\sum_{u\in V}f(u)z_u=\sum_{u\in B}c(u)z_u$ vanishes outside $B$ by Step~1 and vanishes on $B$ by assumption; hence $\sum_u f(u)z_u=0$. By Lemma~\ref{l6.7}, the function $g$ constructed from $f$ is constant on the component $B$: $g\equiv t$ on $B$. This means that $c(v)=t$ for $v\in L\cap B$ and $c(w)=-(k_0/\delta)\,t$ for $w\in U\cap B$, i.e.\ $c=t\,\xi$, where $\xi\colon B\to\mathbb{Q}$ is the fixed vector equal to $1$ on $L\cap B$ and to $-k_0/\delta$ on $U\cap B$. Thus the right kernel of $Z_B$ is contained in $\mathbb{Q}\xi$.
		\par \textbf{Step 3:} the left kernel of $Z_B$ equals $\mathbb{Q}\mathbf{1}_B$. The rank of a matrix equals the rank of its transpose, so the left kernel $\{r\colon B\to\mathbb{Q}:\ \sum_{x\in B}r(x)z_u(x)=0\ \text{for all}\ u\in B\}$ has the same dimension as the right kernel, i.e.\ at most $1$. On the other hand, for $u\in B$, by Step~1,
		$$ \sum_{x\in B}z_u(x)=\sum_{x\in V}z_u(x)=0, $$
		so $\mathbf{1}_B\ne 0$ lies in the left kernel. Consequently, the left kernel equals $\mathbb{Q}\mathbf{1}_B$.
		\par \textbf{Step 4:} $B$ is a module of the graph $G$. Let $i\in V\setminus B$ and $u\in B$. By \ref{eq2}
		$$ (A_Gz_u)(i)=\lambda_u z_u(i)+k_u\varepsilon_u(i)=0, $$
		since $z_u(i)=0$ (Step~1, $i\notin B$) and $i\ne u$. On the other hand, again by Step~1, $(A_Gz_u)(i)=\sum_{x\in V}A_G(i,x)z_u(x)=\sum_{x\in B}A_G(i,x)z_u(x)$. Hence the row $r_i=\bigl(A_G(i,x)\bigr)_{x\in B}$ lies in the left kernel of $Z_B$, and by Step~3 $r_i=t_i\mathbf{1}_B$ for some $t_i\in\mathbb{Q}$. Since the entries of $r_i$ are $0$ or $1$, the vertex $i$ is adjacent either to all vertices of $B$ or to none of them. Since $i\notin B$ is arbitrary, $B$ is a module of $G$.
		\par \textbf{Step 5: $B=V$.} Suppose that $B\ne V$. Take $u\in B$ and $K=\alt(u)$. $G$ and $K$ have the same deck and differ only on pairs $\{u,x\}$ with $z_u(x)\ne 0$, and for such pairs $x\in B$ by Step~1; hence all such pairs are contained in $B$, and $G$ and $K$ coincide on all pairs not contained in $B$. Since $B\subsetneq V$ is a module of $G$, Theorem~\ref{t6.4} gives $\alt(u)=K\cong G$. But $(G,H)$ is counterexample, hence $\alt(u)\cong H\not\cong G$. This contradiction shows that $B=V$, i.e.\ the graph $\mathcal T$ is connected.
	\end{proof}
	
	\section{Absence of an eigenvalue $\theta$}
	All statements are true for every counterexample. By Lemma \ref{l2.5}, we can replace $(G,H)$ with $(\overline{G},\overline{H})$ and conclude that $\theta \leq 0$, and hence $k_1 \leq k_0$.
	\begin{lem}\label{l7.1}
		If $k_1 \leq k_0$ and $v \in S_1(w)$, then $S_0(v) \cap R_1(w) \neq \varnothing$.
	\end{lem}
	\begin{proof}
		Since $v \nsim w$, we have $w \notin S_0(w)$, and by Lemma \ref{l5.1}(b)
		$$	|S_0(v)\cap R_1(w)|+|R_0(v)\cap S_1(w)| = k_0. $$
		If the first intersection is empty, then $k_0 = |R_0(v)\cap S_1(w)| \leq |S_1(w) \setminus \{v\}|=k_1-1$, but this contradicts $k_1 \leq k_0$.
	\end{proof}
	
	\begin{thm}\label{t7.2}
		$\theta$ is not an eigenvalue of $A_G$.
	\end{thm}
	\begin{proof}
		Let $A_Gf=\theta f$. By Lemma \ref{l4.2}(a), $f=0$ on $L$. For every $v \in L$, choose $s(v) \in S_0(v)$, and define
		$$ val(x) =
		\begin{cases} 
			f(s(x)), & x \in L,\\ 
			f(x), & x \in U.
		\end{cases}
		$$
		We show that $f$ is constant on every edge of $\mathcal T$.
		\begin{enumerate}
			\item \textbf{$D_0$}: Take $x \in S_0(v)$. By Lemma \ref{l4.2}(b), $val(x)$ is constant on $S_0(v)=S(\alt(v))$, hence $val(v)=f(x)=val(x)$;
			\item \textbf{$R_1$}: Lemma \ref{l4.2}(c);
			\item \textbf{$D_1$}: Take $v \in S_1(w)$. By Lemma \ref{l7.1}, for $x \in S_0(v) \cap R_1(w)$ we have $val(v)=f(x)=f(w)$;
			\item \textbf{$R_0$}: Take $c \in R_0(v)$. By Lemma \ref{l5.1}(c), there exists $w$ such that $c \in S_1(w)$ and $val(w)=val(c)$. Since $\langle z_v, z_w\rangle \geq z_v(c)z_w(c)=1$, Lemma \ref{l5.1} implies that $w \in S_0(v)$ or $v \in S_1(w)$; in both cases $val(v)=val(w)$.
		\end{enumerate}
		By Theorem \ref{t6.4}, $val(x)\equiv C$, i.e., $f\equiv C$ on $U$ and $f\equiv0$ on $L$. Since the eigenvector is nonzero, there exists $v_0 \in L$ such that
		$$0=\theta f(v_0)=(A_G f)(v_0)=bC,$$
		where $b \geq 1$. Thus $C= 0$ and $f = 0$.
	\end{proof}
	
	\section{$R_1$ connectivity}
	\begin{lem}\label{l8.1}
		Let $g : U \to \mathbb{Q}$ such that $R_1g=k_1g$. Then $D_0^{\mathsf T}D_1g=k_0k_1g$ and $R_0D_1g=k_0D_1g$.
	\end{lem}
	\begin{proof}
		Put $c: V \to \mathbb{Q}$ such that $c_L=-\frac{D_1g}{k_0}$ and $c_U=g$. Denote $W=\sum_{u \in V}c(u)z_u$.
		\par By \ref{eq2} 
		$$
		(A_G-\theta I)W=-\sum_{u \in U}g(u)z_u + \sum_{u \in V}c(u)k_u\varepsilon_u.
		$$
		For $y \in U$:
		$$
		(A_G-\theta I)W(y)=-\sum_{u\in U}g(u)z_u(y)+g(y)k_1 = -R_1g(y)+k_1g(y)=0
		$$
		by lemma conditions.
		For $y \in L$:
		$$
		(A_G-\theta I)W(y)=-\sum_{u\in U}g(u)z_u(y)-D_1g(y).
		$$
		Since $z_u(y)=-\textbf{1}_{S_1(u)}(y)$ we can write $-\sum_{u\in U}g(u)z_u(y)=\sum_{u\in U}g(u)\textbf{1}_{S_1(u)}(y)=D_1g(y)$, hence $(A_G-\theta I)W(y)=0$. 
		\par By Theorem \ref{t7.2} $W\equiv0$ and by Lemma \ref{l6.5} the equations
		$$
		R_0c_L=k_0c_L \qquad D_0^{\mathsf T}c_L=-k_1c_U
		$$
		is exactly what we need.
	\end{proof}
	
	Consider the graph on $U$ with the relation $R_1$, i.e., $w \sim w' \iff w' \in R_1(w)$.
	\begin{thm}\label{t8.2}
		The graph of the relation $R_1$ on $U$ is connected.
	\end{thm}
	\begin{proof}
		Let $C \subseteq U$ and $g = \mathbf{1}_C$. Then, by Lemma \ref{l5.1}(c), put $\varphi_C = \#\{a \in C: u \in S_1(a)\}$. By Lemma \ref{l8.1}, we have
		\begin{enumerate}
			\item[(i)] If $\varphi_C(u) \ge 0$, then $S_0(u) \subseteq C$ and we have $k_0 k_1 g(y) = (D_0^{\mathsf T}\varphi_C)(y) \ge \varphi_C(y) \ge 0$;
			\item[(ii)] $\varphi_C$ is constant on $R_0$-edges because $R_0\varphi_C = k_0\varphi_C$ by Lemma \ref{l6.5}.
		\end{enumerate}
		
		Now, let $C$ be a connected component of $R_1$ on $U$, and let $C' = U \setminus C$. Define 
		$$
		X = \{u \in L : \varphi_C(u) > 0\}, \quad T = C \cup X.
		$$
		We will show that $T$ is closed with respect to the edges of $\mathcal{T}$ (in both directions). The idea is to verify that for any vertex in $T$, all its neighbors in $\mathcal{T}$ also lie in $T$.
		
		\begin{itemize}
			\item[(a)] \textit{Outgoing edges from $a \in C$:} 
			Since $C$ is a component of $R_1$, all $R_1$-neighbors of $a$ remain in $C$, so $R_1(a) \subseteq C$. Moreover, for any $u \in S_1(a)$, the vertex $u$ is adjacent to $a \in C$, which implies $\varphi_C(u) \ge 1 > 0$. Thus, $S_1(a) \subseteq X$.
			
			\item[(b)] \textit{Outgoing edges from $u \in X$:} 
			By property (i), since $\varphi_C(u) > 0$, we have $S_0(u) \subseteq C$. By property (ii), the function $\varphi_C$ is constant along $R_0$-edges, so $R_0(u) \subseteq X$.
			
			\item[(c)] \textit{Incoming edges to $a \in C$:} 
			Consider an edge in $\mathcal{T}$ incident to $a \in C$, and let its other endpoint be $c$. We show $c \in T$.
			\begin{itemize}
				\item[(1)] If $c \in U$ and $a \in R_1(c)$, then by symmetry $c \in R_1(a) \subseteq C$.
				\item[(2)] If $c \in L$ and $a \in S_0(c)$ (which means $c \in S_1(a)$), suppose for contradiction that $c \notin X$. Then $\varphi_C(c) = 0$. By Lemma \ref{l5.1}(c), there exists $a' \in C'$ such that $c \in S_1(a')$. This implies $\varphi_{C'}(c) > 0$. By property (i) applied to $C'$, we get $S_0(c) \subseteq C'$. But this contradicts $a \in S_0(c) \cap C$. Therefore, $c \in X \subseteq T$.
			\end{itemize}
			
			\item[(d)] \textit{Incoming edges to $u \in X$:} 
			Consider an edge in $\mathcal{T}$ incident to $u \in X$, and let its other endpoint be $c$.
			\begin{itemize}
				\item[(1)] If $c \in L$ and $u \in R_0(c)$, then by symmetry $c \in R_0(u)$. By property (ii), $c \in X$.
				\item[(2)] If $c \in U$ and $u \in S_1(c)$, suppose for contradiction that $c \notin C$ (so $c \in C'$). Since $u \in S_1(c)$, we have $\varphi_{C'}(u) > 0$. However, $u \in X$ also implies $\varphi_C(u) > 0$. Applying property (i) to both $C$ and $C'$, we obtain $S_0(u) \subseteq C$ and $S_0(u) \subseteq C'$. This means $S_0(u) \subseteq C \cap C' = \varnothing$, which is impossible because $k_0 \ge 1$ ensures $S_0(u)$ is non-empty. Thus, $c \in C \subseteq T$.
			\end{itemize}
		\end{itemize}
		
		We have shown that $T$ is a non-empty subset of vertices that is completely closed under the edges of $\mathcal{T}$. By Theorem \ref{t6.8}, the graph $\mathcal{T}$ is connected, which forces $T = V$. Since $T = C \cup X$ and $X \subseteq L$, it follows that $C$ must coincide with all of $U$. Therefore, the graph of the relation $R_1$ on $U$ is connected.
	\end{proof}
	
	\begin{thm}\label{t8.3}
		For every connected component $L_i$ of the graph $R_0$ and every vertex $w \in U$, there exists a vertex $v' \in L_i$ such that $w \in S_0(v')$.
	\end{thm}
	\begin{proof}
		Fix a connected component $L_i$ of the graph $R_0$. We begin by defining a counting function for the number of neighbors a vertex $y \in U$ has in $L_i$ via the $S_0$-relation. Let
		$$
		\varphi(y) = \#\{u \in L_i : y \in S_0(u)\} = (D_0^{\mathsf T} \mathbf{1}_{L_i})(y), \quad y \in U.
		$$
		Our goal is to show that $\varphi(w) > 0$ for every $w \in U$.
		
		To do this, we take the identity from Lemma \ref{l5.1}(c) and sum it over all $u \in L_i$ for a fixed vertex $w' \in U$. This yields:
		$$
		\sum_{u \in L_i} |S_0(u) \cap R_1(w')| + \sum_{u \in L_i} |R_0(u) \cap S_1(w')| = k_1 \sum_{u \in L_i} [w' \in S_0(u)] + k_0 \sum_{u \in L_i} [u \in S_1(w')].
		$$
		
		Let us rewrite the left-hand side by swapping the order of summation. 
		For the first term, summing the intersections over $u \in L_i$ is equivalent to summing over all $x \in R_1(w')$ the number of $u \in L_i$ such that $x \in S_0(u)$. By definition, this is exactly $\sum_{x \in R_1(w')} \varphi(x)$.
		For the second term, we sum over $x \in S_1(w')$ the number of $u \in L_i$ such that $x \in R_0(u)$. Thus, the left-hand side becomes:
		$$
		\sum_{x \in R_1(w')} \varphi(x) + \sum_{x \in S_1(w')} \#\{u \in L_i : x \in R_0(u)\}.
		$$
		
		Now, notice that the right-hand side of our summed equation simplifies nicely. The first sum on the right is simply $k_1 \varphi(w')$, and the second is $k_0 |L_i \cap S_1(w')|$. So we have:
		$$
		\sum_{x \in R_1(w')} \varphi(x) + \sum_{x \in S_1(w')} \#\{u \in L_i : x \in R_0(u)\} = k_1 \varphi(w') + k_0 |L_i \cap S_1(w')|.
		$$
		
		Next, we examine the second sum on the left. Since the relation $R_0$ is symmetric, the condition $x \in R_0(u)$ is equivalent to $u \in R_0(x)$. Therefore, we are counting the $R_0$-neighbors of $x$ that lie in $L_i$. Because $L_i$ is a full connected component of $R_0$, a vertex $x$ has either all $k_0$ of its $R_0$-neighbors in $L_i$ (if $x \in L_i$), or none of them (if $x \notin L_i$). This means:
		$$
		\#\{u \in L_i : x \in R_0(u)\} = |R_0(x) \cap L_i| = k_0 [x \in L_i].
		$$
		
		Substituting this back into our sum, the second term on the left becomes \\
		$\sum_{x \in S_1(w')} k_0 [x \in L_i]$, which is exactly $k_0 |S_1(w') \cap L_i|$. 
		
		Notice that this perfectly matches the second term on the right-hand side. They cancel each other out, leaving us with:
		$$
		\sum_{x \in R_1(w')} \varphi(x) = k_1 \varphi(w').
		$$
		In matrix notation, this is written as $R_1 \varphi = k_1 \varphi$. This equation tells us that $\varphi$ is an eigenfunction of the adjacency matrix $R_1$ corresponding to the eigenvalue $k_1$ (which is the degree of the graph). 
		
		A fundamental property of any connected $k_1$-regular graph is that the eigenspace for the eigenvalue $k_1$ is one-dimensional and consists only of constant vectors. Since we proved in Theorem \ref{t8.2} that the graph $R_1$ on $U$ is connected, $\varphi$ must be constant on all of $U$.
		
		Finally, we note that the sum of $\varphi$ over all $y \in U$ is strictly positive:
		$$
		\sum_{y \in U} \varphi(y) = \sum_{y \in U} \#\{u \in L_i : y \in S_0(u)\} = \sum_{u \in L_i} |S_0(u)| = k_0 |L_i| > 0.
		$$
		Because the sum is positive and $\varphi$ is constant on $U$, the constant value must be strictly positive. Thus, $\varphi(w) > 0$ for all $w \in U$. 
		
		By the definition of $\varphi$, the fact that $\varphi(w) > 0$ means that there is at least one vertex $v' \in L_i$ such that $w \in S_0(v')$.
	\end{proof}
	
	\section{Contradiction with counterexample}
	
	\begin{thm}\label{t9.1}
		$D_2 = D_3 = 0$.
	\end{thm}
	
	\begin{proof}
		Let $\Omega_i = \operatorname{supp} D_i \subseteq L \times U$. Since all matrices are non-negative, any positive summand makes the corresponding entry of the matrix product positive. Throughout, we assume $v' \in R_0(v)$ or $w' \in R_1(w)$:
		\begin{itemize}
			\item A step from $\Omega_1$ along $R_0$ lands in $\Omega_0 \cup \Omega_1$. Indeed, if $(v, w) \in \Omega_1$, then $(R_0D_1)(v', w) \ge 1$, and Lemma \ref{l5.1}(b) at $(v', w)$ implies $D_0(v', w) + D_1(v', w) > 0$.
			
			\item A step from $\Omega_0$ along $R_0$ cannot land in $\Omega_3$. Otherwise, the left-hand side of Proposition \ref{p5.3}(a) at $(v, w)$ would be positive, while the right-hand side is zero.
			
			\item A step from $\Omega_0$ along $R_1$ lands in $\Omega_0 \cup \Omega_1$, which follows directly from $(D_0R_1)(v, w') \ge 1$ and Lemma \ref{l5.1}(b).
			
			\item A step from $\Omega_2$ along $R_1$ cannot land in $\Omega_0$. If it did, we would have $(D_0R_1)(v, w) \ge 1$, contradicting the right-hand side of Lemma \ref{l5.1}(b) being zero at $(v, w) \in \Omega_2$.
			
			\item A step from $\Omega_2$ along $R_1$ cannot land in $\Omega_1$. Similarly, this would yield $(D_2R_1)(v, w') \ge 1$, contradicting the right-hand side of Proposition \ref{p5.3}(b) being zero at $(v, w') \in \Omega_1$.
		\end{itemize}
		
		Now we claim that a step from $\Omega_0$ along $R_0$ leads to $\Omega_0 \cup \Omega_1$. 
		
		Let $(v, w) \in \Omega_0$ and $v' \in R_0(v)$, and suppose for contradiction that $(v', w) \notin \Omega_0 \cup \Omega_1$. By the second observation above, this implies $(v', w) \in \Omega_2$. 
		
		The row $v'$ of $D_0$ is non-zero; thus, there exists a vertex $w_0$ such that $(v', w_0) \in \Omega_0$. By Theorem \ref{t8.2}, there exists a path $w = x_0, x_1, \dots, x_m = w_0$ in $R_1$. We prove by induction on $i$ that $(v, x_i) \in \Omega_0$ and $(v', x_i) \in \Omega_2$. 
		
		Assume this holds for some $i$. By the fourth and fifth observations, $(v', x_{i+1}) \notin \Omega_0 \cup \Omega_1$. By the third observation, $(v, x_{i+1}) \in \Omega_0 \cup \Omega_1$. The case $(v, x_{i+1}) \in \Omega_1$ is impossible by the first observation; hence, $(v, x_{i+1}) \in \Omega_0$. Furthermore, $(v', x_{i+1}) \notin \Omega_3$ by the second observation, which forces $(v', x_{i+1}) \in \Omega_2$. 
		
		For $i = m$, this gives $(v', w_0) \in \Omega_2$, contradicting the choice of $(v', w_0) \in \Omega_0$. Thus, $\Omega_0 \cup \Omega_1$ is closed under steps along $R_0$.
		
		Finally, by Theorem \ref{t8.3}, for any $v, w$, in the $R_0$-component containing $v$, there exists $v'$ such that $(v', w) \in \Omega_0$. Moving from $v'$ to $v$ along $R_0$ keeps us within the closed set $\Omega_0 \cup \Omega_1$, which yields $(v, w) \in \Omega_0 \cup \Omega_1$. Hence, $\Omega_2 = \Omega_3 = \varnothing$.
	\end{proof}
	
	\begin{lem}\label{l9.2}
		Let $F$ be an $m$-regular graph and let $r$ be a vertex of $F$. If all cards $F-v$ for $v \neq r$ are pairwise isomorphic, then the card $F-r$ is also isomorphic to them.
	\end{lem}
	
	\begin{proof}
		Suppose, for the sake of contradiction, that all cards $F-v$ for $v \neq r$ are pairwise isomorphic, but $F-r$ is not isomorphic to them.
		
		In an $m$-regular graph, the neighbors of a deleted vertex $v$ are exactly the vertices of degree $m-1$ in the card $F-v$. Therefore, any isomorphism $F-v \cong F-w$ maps the neighborhood of $v$ to the neighborhood of $w$ and extends (via $v \mapsto w$) to an automorphism of $F$. 
		
		Thus, for any $v, w \neq r$, there exists an automorphism $\pi \in \operatorname{Aut}(F)$ such that $\pi(v) = w$. If $\pi(r) \neq r$, then $\pi$ induces an isomorphism $F-r \cong F-\pi(r)$, which contradicts our assumption that $F-r$ is unique. Therefore, $\pi(r) = r$, and since $\pi$ preserves adjacency, we get $r \sim v \iff r \sim w$. This means that $r$ is either adjacent to all other vertices in $F$ or to none of them.
		
		Both of these cases are impossible. If $r$ is adjacent to all other vertices in $F$, then $F \cong K_{m+1}$, but all cards in the deck of a complete graph are pairwise isomorphic. If $r$ is an isolated vertex, then $F$ is a graph on $m+1$ vertices without edges, and all cards in its deck are pairwise isomorphic.
	\end{proof}
	
		\begin{thm}
		Equality $D_2=D_3=0$ are impossible
	\end{thm}
	
	\begin{proof}
		From $D_2=0$, we obtain $S_0(v)=N_G(v) \cap U$ for all $v \in L$, i.e., $k_0=b$.
		\par From $D_3=0$, we obtain $S_1(w)=L \setminus N_G(w)$ for all $w \in U$, i.e., $k_1 = p - c$. By Lemma \ref{l5.2},
		$$
		k_1 - k_0 = \theta = p - c - b = e - b,
		$$
		thus $p=c+e=d+1$. For $v \in L$, the number of low non-neighbors is equal to $p-1-(d-b)=b=k_0$, hence $R_0(v) = L \setminus (N_G(v) \cup \{v\})$. For $w \in U$, we have $|N_G(w) \cap U|=e=p-c=k_1$, thus $R_1(w)=N_G(w) \cap U$.
		\par Now we show that in $\alt(u)$, the vertex $u$ is adjacent to all vertices in $L \setminus \{u\}$ if $u \in L$, and to all vertices in $L$ if $u \in U$.
		\begin{enumerate}
			\item Let $u \in L$. In $G$, the vertex $u$ is adjacent to the vertices in
			$$
			(N_G(u) \cap L) \cup (N_G(u) \cap U) =(N_G(u) \cap L) \cup S_0(u).
			$$
			Then, in $\alt(u)$, the vertex $u$ is adjacent to the vertices in
			$$
			(N_G(u) \cap L) \cup R_0(u)=(N_G(u) \cap L) \cup (L \setminus (N_G(u) \cup \{u\}))=L \setminus \{u\}.
			$$
			\item Let $u \in U$. In $G$, the vertex $u$ is adjacent to the vertices in
			$$
			(N_G(u) \cap L) \cup (N_G(u) \cap U) =(N_G(u) \cap L) \cup R_1(u).
			$$
			Then, in $\alt(u)$, the vertex $u$ is adjacent to the vertices in
			$$
			(N_G(u) \cap L) \cup S_1(u)=(N_G(u) \cap L) \cup (L \setminus N_G(u))=L.
			$$
		\end{enumerate}
		Add a new vertex $r$ to $G$, adjacent only to the vertices in $L$. We construct a new graph $F$, which is $(d+1)$-regular, since $\deg r = p = d+1$. We have $F - r \cong G$. For $v \in V$,
		$$
		G-v \cong (F-r)-v \cong (F-v)-r,
		$$
		hence $F-v$ is one of the extensions and, apart from the vertex $r$, coincides with $G$. Now we show that $F-v \cong \alt(r)$ for all $v \in V$.
		\par Indeed, if $v \in L$, then $r$ is low in $F-v$ and adjacent to all vertices in $L \setminus \{r\}$. If $v \in U$, then $r$ is high and adjacent to all vertices in $L$. Hence, as shown above, $F - v \cong \alt(t) \cong H$ for all $v \in V$. By Lemma \ref{l9.2}, $F - v \cong F -r$, but $G \ncong H$. This is a contradiction, so no counterexample exists.
	\end{proof}

\end{document}